\documentclass[12pt]{amsart}
\usepackage{latexsym, amsmath, amscd, amssymb, amsthm}
\usepackage{bm, mathrsfs}

\usepackage{calc}
\DeclareMathOperator{\BBox}{Box}
\DeclareMathOperator{\Coker}{Coker}
\DeclareMathOperator{\cone}{Cone}
\DeclareMathOperator{\dg}{DG}
\DeclareMathOperator{\Ext}{Ext}
\DeclareMathOperator{\ev}{ev}
\DeclareMathOperator{\Hilb}{Hilb}
\DeclareMathOperator{\Hom}{Hom}
\DeclareMathOperator{\image}{Im}
\DeclareMathOperator{\Ker}{Ker}
\DeclareMathOperator{\link}{link}
\DeclareMathOperator{\pr}{pr}
\DeclareMathOperator{\rank}{rank} 
\DeclareMathOperator{\Spec}{Spec}
\DeclareMathOperator{\Sym}{Sym}
\DeclareMathOperator{\Pic}{Pic}

\newtheorem{lemma}{Lemma}[section]
\newtheorem{theorem}[lemma]{Theorem}
\newtheorem{corollary}[lemma]{Corollary}
\newtheorem{proposition}[lemma]{Proposition}
\theoremstyle{definition}

\newtheorem{condition}[lemma]{Condition}
\newtheorem{example}[lemma]{Example}

\newtheorem{remark}[lemma]{Remark}
\newtheorem*{convention}{Conventions}
\newtheorem*{acknowledgment}{Acknowledgments}
\theoremstyle{remark}
\renewcommand{\theequation}%
{\arabic{section}.\arabic{lemma}.\arabic{equation}}
\begin{document}

\title[The orbifold {C}how ring of toric {D}eligne-{M}umford
stacks]{The orbifold {C}how ring of \\ toric {D}eligne-{M}umford
  stacks}

\author{Lev A.  Borisov}
\address{Department of Mathematics \\ University of Wisconsin \\
  Madison \\ WI 53706 \\ USA}
\email{borisov@math.wisc.edu}

\author{Linda Chen}
\address{Department of Mathematics \\ Columbia University \\
  New York \\ NY 10027 \\ USA} 
\email{lchen@math.columbia.edu}

\author{Gregory G. Smith}
\address{Department of Mathematics \\ Barnard College \\
  Columbia University \\ New York \\ NY 10027 \\ USA}
\email{ggsmith@math.columbia.edu}  

\thanks{The first author was partially supported by NSF grant
  DMS-0140172 and the second author was partially supported by NSF
  VIGRE grant DMS-9810750}

\subjclass[2000]{Primary 14N35; Secondary 14C15, 14M25}

\begin{abstract}
  Generalizing toric varieties, we introduce toric Deligne-Mumford
  stacks.  The main result in this paper is an explicit calculation of
  the orbifold Chow ring of a toric Deligne-Mumford stack.  As an
  application, we prove that the orbifold Chow ring of the toric
  Deligne-Mumford stack associated to a simplicial toric variety is a
  flat deformation of (but is not necessarily isomorphic to) the Chow
  ring of a crepant resolution.
\end{abstract}

\maketitle

\vspace*{-1em}

\section{Introduction} \label{sec:intro}

The orbifold Chow ring of a Deligne-Mumford stack, defined by
Abramovich, Graber and Vistoli~\cite{AGV}, is the algebraic version of
the orbifold cohomology ring introduced by W.~Chen and
Ruan~\cite{ChenRuan0}~\cite{ChenRuan}.  By design, this ring
incorporates numerical invariants, such as the orbifold Euler
characteristic and the orbifold Hodge numbers, of the underlying
variety.  The product structure is induced by the degree zero part of
the quantum product; in particular, it involves Gromov-Witten
invariants.  Inspired by string theory and results in
Batyrev~\cite{Batyrev} and Yasuda~\cite{Yasuda}, one expects that, in
nice situations, the orbifold Chow ring coincides with the Chow ring
of a resolution of singularities.  Fantechi and
G\"{o}ttsche~\cite{FantechiGottsche} and Uribe~\cite{Uribe} verify
this conjecture when the orbifold is $\Sym^{n}(S)$ where $S$ is a
smooth projective surface with $K_{S} = 0$ and the resolution is
$\Hilb^{n}(S)$.  The initial motivation for this project was to
compare the orbifold Chow ring of a simplicial toric variety with the
Chow ring of a crepant resolution.

To achieve this goal, we first develop the theory of toric
Deligne-Mumford stacks.  Modeled on simplicial toric varieties, a
toric Deligne-Mumford stack corresponds to a combinatorial object
called a \emph{stacky fan}.  As a first approximation, this object is
a simplicial fan with a distinguished lattice point on each ray in the
fan.  More precisely, a stacky fan $\bm{\Sigma}$ is a triple
consisting of a finitely generated abelian group $N$, a simplicial fan
$\Sigma$ in $\mathbb{Q} \otimes_{\mathbb{Z}} N$ with $n$ rays, and a
map $\beta \colon \mathbb{Z}^{n} \rightarrow N$ where the image of the
standard basis in $\mathbb{Z}^{n}$ generates the rays in $\Sigma$.  A
rational simplicial fan $\Sigma$ produces a canonical stacky fan
$\bm{\Sigma} := (N,\Sigma,\beta)$ where $N$ is the distinguished
lattice and $\beta$ is the map defined by the minimal lattice points
on the rays.  Hence, there is a natural toric Deligne-Mumford stack
associated to every simplicial toric variety.  A stacky fan
$\bm{\Sigma}$ encodes a group action on a quasi-affine variety and the
toric Deligne-Mumford stack $\mathcal{X}(\bm{\Sigma})$ is the
quotient.  If $\Sigma$ corresponds to a smooth toric variety
$X(\Sigma)$ and $\bm{\Sigma}$ is the canonical stacky fan associated
to $\Sigma$, then we simply have $\mathcal{X}(\bm{\Sigma}) =
X(\Sigma)$.  We show that many of the basic concepts, such as open and
closed toric substacks, line bundles, and maps between toric
Deligne-Mumford stacks, correspond to combinatorial notions.  We
expect that many more results about toric varieties to lift to the
realm of stacks and hope that toric Deligne-Mumford stacks will serve
as a useful testing ground for general theories.

Our description for the orbifold Chow ring of a toric Deligne-Mumford
stack $\mathcal{X}(\bm{\Sigma})$ parallels the ``Stanley-Reisner''
presentation for the Chow ring of a simplicial toric variety.
Specifically, the stacky fan $\bm{\Sigma}$ gives rise to the
\emph{deformed group ring} $\mathbb{Q}[N]^{\bm{\Sigma}}$.  As a
$\mathbb{Q}$-vector space, $\mathbb{Q}[N]^{\bm{\Sigma}}$ is the group
algebra of $N$.  Since $N$ is abelian, we write
$\mathbb{Q}[N]^{\bm{\Sigma}} = \bigoplus_{c \in N} \mathbb{Q} \,
y^{c}$ where $y$ is a formal variable.  For $c \in N$, $\bar{c}$
denotes the image of $c$ in $\mathbb{Q} \otimes_{\mathbb{Z}} N$.
Multiplication in $\mathbb{Q}[N]^{\bm{\Sigma}}$ is defined by the
equation:
\[
y^{c_{1}} \cdot y^{c_{2}} :=
\begin{cases}
  y^{c_{1} + c_{2}} & \text{ if there is $\sigma \in \Sigma$ such that
    $\bar{c}_{1} \in \sigma$ and $\bar{c}_{2} \in \sigma$,}\\
  0 & \text{otherwise.}
\end{cases}
\]
Let $b_{i}$ be the image under the map $\beta$ of the $i$-th standard
basis vector.  The map $\beta \colon \mathbb{Z}^{n} \rightarrow N$
endows $\mathbb{Q}[N]^{\bm{\Sigma}}$ with a $\mathbb{Q}$-grading; if
$\bar{c} = \sum_{\bar{b}_{i} \in \sigma} m_{i} \bar{b}_{i}$ where
$\sigma$ is the minimal cone in $\Sigma$ containing $\bar{c}$ and
$m_{i}$ is a nonnegative rational number, then the
$\mathbb{Q}$-grading is given by $\deg(y^{c}) := \sum_{\bar{b}_{i} \in
  \sigma} m_{i}$.  Writing $A_{orb}^{*} \bigl(
\mathcal{X}(\bm{\Sigma}) \bigr)$ for the orbifold Chow ring of
$\mathcal{X}(\bm{\Sigma})$ with rational coefficient, our principal
result is:

\begin{theorem} \label{thm:main}
  If $\mathcal{X}(\bm{\Sigma})$ is a complete toric Deligne-Mumford
  stack, then there is an isomorphism of $\mathbb{Q}$-graded rings:
  \[
  A_{orb}^{*} \bigl( \mathcal{X}(\bm{\Sigma}) \bigr) \cong
  \frac{\mathbb{Q}[N]^{\bm{\Sigma}}}{\bigl\langle \sum_{i=1}^{n}
    \theta(b_{i}) y^{b_{i}} : \theta \in \Hom(N, \mathbb{Z})
    \bigr\rangle} \, .
  \]
\end{theorem}

\noindent
Using differential geometry, Jiang~\cite{Jiang} establishes this
result for the weighted projective space $\mathbb{P}(1,2,2,3,3,3)$.

Our proof of this theorem involves two steps.  By definition, the
orbifold Chow ring $A_{orb}^{*} \bigl( \mathcal{X}(\bm{\Sigma})
\bigr)$ is isomorphic as an abelian group to the Chow ring of the
inertia stack $\mathcal{I} \bigl( \mathcal{X}(\bm{\Sigma}) \bigr)$.
We first express $\mathcal{I} \bigl( \mathcal{X}(\bm{\Sigma}) \bigr)$
as a disjoint union of toric Deligne-Mumford stacks.  This leads to a
proof of Theorem~\ref{thm:main} at the level of $\mathbb{Q}$-graded
vector spaces.  To compare the ring structures, we also express the
moduli space $\mathcal{K}_{0,3} \bigl( \mathcal{X}(\bm{\Sigma}), 0
\bigr)$ of $3$-pointed twisted stable maps as a disjoint union of
toric Deligne-Mumford stacks.  This combinatorial description allows
us to compute the virtual fundamental class of $\mathcal{K}_{0,3}
\bigl( \mathcal{X}(\bm{\Sigma}), 0 \bigr)$.  We are then able to
verify that multiplication in the deformed group ring coincides with
the product in the orbifold Chow ring.

The paper is organized as follows.  In Section~\ref{sec:gale}, we
extend Gale duality to maps of finitely generated abelian groups.
This duality forms an essential link between stacky fans and toric
Deligne-Mumford stacks.  Nevertheless, this theory is entirely
self-contained and may be of interest in other situations.  The
rudimentary theory of toric Deligne-Mumford stacks is developed in
Sections~\ref{sec:stackytoric} and~\ref{sec:openclosed}.
Specifically, we detail the correspondence between stacky fans and
toric Deligne-Mumford stacks, we describe the open and closed toric
substacks and we express the inertia stacks as disjoint unions of
toric Deligne-Mumford stacks.  The proof of Theorem~\ref{thm:main} is
given in Sections~\ref{sec:chowanddef} and~\ref{sec:orbchow}.  Finally
in Section~\ref{sec:crepant}, we use our main result to compare the
orbifold Chow rings of a simplicial toric variety and its crepant
resolutions.

\begin{convention}
  Throughout this paper, we work over the field $\mathbb{C}$ of
  complex numbers and consider Chow rings and orbifold Chow rings with
  rational coefficients.
\end{convention}

\begin{acknowledgment}
  We would like to thank Dan Abramovich, Bob Friedman, Bill Fulton,
  Tom Graber, Paul Horja, Andrew Kresch, Martin Olsson, Michael
  Thaddeus, Hsian-Hua Tseng, Howard Thompson and Ravi Vakil for useful
  discussions.
\end{acknowledgment}

\section{Gale Duality with Torsion} \label{sec:gale}
\setcounter{equation}{0}

In this section, we generalize Gale duality.  To orient the reader, we
first recall the basic form of Gale duality and its application to
toric geometry.  Given $n$ vectors $b_{1}, \dotsc, b_{n}$ which span
$\mathbb{Q}^{d}$, there is a dual configuration $[a_{1} \dotsb a_{n}]
\in \mathbb{Q}^{(n-d) \times n}$ such that
\begin{equation} \label{equ:dualses}
0 \longrightarrow \mathbb{Q}^{d} \xrightarrow{[ b_{1} \dotsb
  b_{n}]^{\textsf{T}}} \mathbb{Q}^{n} \xrightarrow{[a_{1} \dotsb
  a_{n}]} \mathbb{Q}^{n-d} \longrightarrow 0
\end{equation}
is a short exact sequence; see Theorem~6.14 in \cite{Ziegler}. The set
of vectors $\{ a_{1}, \dotsc, a_{n} \}$ is uniquely determined up to a
linear coordinates transformation in $\mathbb{Q}^{n-d}$.  This duality
plays a role in study of smooth toric varieties.  Specifically, let
$\Sigma$ be a fan with $n$ rays such that the corresponding toric
variety $X(\Sigma)$ is smooth.  If $N \cong \mathbb{Z}^{d}$ is the
lattice in $\Sigma$, then the minimal lattice points $b_{1}, \dotsc,
b_{n}$ generating the rays determine a map $\beta \colon
\mathbb{Z}^{n} \rightarrow N$.  By tensoring with $\mathbb{Q}$, we
obtain a dual configuration $\{ a_{1}, \dotsc, a_{n} \}$.  Since
$X(\Sigma)$ is smooth, we have $a_{i} \in \mathbb{Z}^{n-d}$ and the
set $\{ a_{1}, \dotsc, a_{n} \}$ is unique up to unimodular
(determinant $\pm 1$) coordinate transformations of
$\mathbb{Z}^{n-d}$.  Abbreviating $\Hom_{\mathbb{Z}}(-,\mathbb{Z})$ by
$(-)^{\star}$, it follows that the set $\{ a_{1}, \dotsc, a_{n} \}$
defines a map $\beta^{\vee} \colon (\mathbb{Z}^{n})^{\star}
\rightarrow\mathbb{Z}^{n-d} \cong \Pic(X)$ and the short exact
sequence \eqref{equ:dualses} becomes $0 \rightarrow N^{\star}
\xrightarrow{\;\; \beta^{\star} \;\;} (\mathbb{Z}^{n})^{\star}
\xrightarrow{\;\; \beta^{\vee} \;\;} \Pic(X) \rightarrow 0$; see
Section~3.4 in \cite{Fulton}.  Our goal is to extend this theory to a
larger class of maps.

Let $N$ be a finitely generated abelian group and consider a group
homomorphism $\beta \colon \mathbb{Z}^{n} \rightarrow N$.  The map
$\beta$ is determined by a finite subset $\{ b_{1}, \dotsc, b_{n} \}$
of $N$.  Throughout this paper, we simply write $(-)^{\star}$ for the
functor $\Hom_{\mathbb{Z}}(-,\mathbb{Z})$.  The dual map $\beta^{\vee}
\colon (\mathbb{Z}^{n})^{\star} \rightarrow \dg(\beta)$ is defined as
follows.  Choose projective resolutions $\bm{E}$ and $\bm{F}$ of the
$\mathbb{Z}$-modules $\mathbb{Z}^{n}$ and $N$ respectively.
Theorem~2.2.6 in \cite{Weibel} shows that $\beta \colon \mathbb{Z}^{n}
\rightarrow N$ lifts to a morphism $\bm{E} \rightarrow \bm{F}$ and
Subsection~1.5.8 in \cite{Weibel} shows that the mapping cone
$\cone(\beta)$ fits into an exact sequence of cochain complexes $0
\rightarrow \bm{F} \rightarrow \cone(\beta) \rightarrow \bm{E}[1]
\rightarrow 0$.  Since $\bm{E}$ is projective, we have the exact
sequence of cochain complexes
\begin{equation} \label{eq:cochainses}
0 \longrightarrow \bm{E}[1]^{\star} \longrightarrow
\cone(\beta)^{\star} \longrightarrow \bm{F}^{\star} \longrightarrow 0
\end{equation}
and the associated long exact sequence in cohomology contains the
exact sequence:
\begin{equation} \label{eq:dualmap}
N^{\star} \xrightarrow{\;\; \beta^{\star} \;\;}
(\mathbb{Z}^{n})^{\star} \longrightarrow H^{1}\big(
\cone(\beta)^{\star} \big) \longrightarrow \Ext_{\mathbb{Z}}^{1}(N,
{\mathbb{Z}}) \longrightarrow 0 \, .
\end{equation}
Set $\dg(\beta) := H^{1}\big( \cone(\beta)^{\star} \big)$ and define
the dual map $\beta^{\vee} \colon (\mathbb{Z}^{n})^{\star} \rightarrow
\dg(\beta)$ to be the second map in \eqref{eq:dualmap}.  Since
$\mathbb{Z}^{n}$ is projective, $\beta^{\vee}$ is in fact the only
nontrivial connecting homomorphism in the long exact sequence
associated \eqref{eq:cochainses}.  This abstract definition guarantees
the naturality of this construction.  Indeed, mapping cones are
natural in the following sense: for every commutative diagram of
cochain complexes:
\[
\begin{CD}
  \bm{E} @>{\beta}>> \bm{F} & \\
  @VV{\nu}V @VV{\upsilon}V & \\
  \bm{E}' @>{\beta'}>> \bm{F}' & \, ,
\end{CD} 
\]
the map $\cone(\beta) \rightarrow \cone(\beta')$ given by $(b,c)
\mapsto \bigl( \nu(b), \upsilon(c) \bigr)$ is a morphism and it is an
isomorphism if $\nu$ and $\upsilon$ are.  It follows that both
$\dg(\beta)$ and $\beta^{\vee}$ are well-defined up to natural
isomorphism. 

On the other hand, there is an explicit description of the dual map
$\beta^{\vee}$ and the dual group $\dg(\beta)$.  If $N$ has rank $d$,
then the structure theorem of finitely generated abelian groups
implies that there exists an integer matrix $Q$ such that $0
\rightarrow \mathbb{Z}^{r} \xrightarrow{Q} \mathbb{Z}^{d+r}
\rightarrow 0$ is a projective resolution of $N$.  The map $\beta
\colon \mathbb{Z}^{n} \rightarrow N$ lifts to a map $\mathbb{Z}^{n}
\rightarrow \mathbb{Z}^{d+r}$ given by a matrix $B$.  Since
$\mathbb{Z}^{n}$ is projective, the cochain complex with $E^{0} =
\mathbb{Z}^{n}$ and $E^{i} = 0$ for all $i \neq 0$ is a projective
resolution of $\mathbb{Z}^{n}$.  With these choices, $\cone(\beta)$ is
the complex $0 \rightarrow \mathbb{Z}^{n + r} \xrightarrow{[ B \; Q ]}
\mathbb{Z}^{d + r} \rightarrow 0$ and we obtain the sequence
\eqref{eq:dualmap} by applying the Snake Lemma to the diagram:
\[
\begin{CD}
  && 0 @>>> (\mathbb{Z}^{d+r})^{\star} @>>>
  ({\mathbb{Z}}^{d+r})^{\star} @>>> 0\\
  && @VVV @VV{[B\;Q]^{\star}}V @VV{Q^{\star}}V \\
  0 @>>> (\mathbb{Z}^{n})^{\star} @>>> (\mathbb{Z}^{n+r})^{\star} @>>>
  (\mathbb{Z}^{r})^{\star} @>>> 0
\end{CD} \quad .
\]
It follows that $\dg(\beta) = (\mathbb{Z}^{n+r})^{\star} / \image( [ B
\; Q ]^{\star})$ and that the map $\beta^{\vee}$ is the composition of
the inclusion map $(\mathbb{Z}^{n})^{\star} \rightarrow
(\mathbb{Z}^{n+r})^{\star}$ and the quotient map
$(\mathbb{Z}^{n+r})^{\star} \rightarrow \dg(\beta)$.

\begin{example} \label{exa:dual}
  The set $\{ (2, 1), (-3, 0) \} \in \mathbb{Z} \oplus \mathbb{Z}/2
  \mathbb{Z}$ yields a map $\beta \colon \mathbb{Z}^{2} \rightarrow
  \mathbb{Z} \oplus \mathbb{Z}/2 \mathbb{Z}$.  In this case, $Q =
  \left[
    \begin{smallmatrix} 0 \\ 2 \end{smallmatrix}
  \right]$ and $B = \left[
    \begin{smallmatrix} 2 &-3 \\ 1 & 0 \end{smallmatrix}
  \right]$.  Since the vector $\left[
    \begin{smallmatrix} 6 & 4 & -3 \end{smallmatrix}
  \right]^{\star}$ spans the integer kernel of matrix $\left[
    \begin{smallmatrix} 2 & -3 & 0 \\ 1 & 0 & 2 \end{smallmatrix}
  \right]$, we have $\dg(\beta) \cong \mathbb{Z}^{3} / \image \bigl(
  \left[
    \begin{smallmatrix} 2 & -3 & 0 \\ 1 & 0 & 2 \end{smallmatrix}
  \right]^{\star} \bigr) \cong \mathbb{Z}$ and $\beta^{\vee} \colon
  \mathbb{Z}^{2} \rightarrow \mathbb{Z}$ is given by the matrix
  $\left[
    \begin{smallmatrix} 6 & 4 \end{smallmatrix}
  \right]$.
\end{example}

We are especially interested in the map $\beta \colon \mathbb{Z}^{n}
\rightarrow N$ when it has a finite cokernel.  

\begin{proposition} 
  Let $\beta \colon \mathbb{Z}^{n} \rightarrow N$ be a homomorphism of
  finitely generated abelian groups.  The map $\beta$ is naturally
  isomorphic to $\beta^{\vee \vee}$ if and only if the cokernel of
  $\beta$ is finite.  Moreover, if cokernel of $\beta$ is finite, then
  the kernel of $\beta^{\vee}$ is $N^{\star}$.
\end{proposition}

\begin{proof}
  Suppose that $\Coker(\beta)$ is not finite.  The sequence
  \eqref{eq:dualmap} implies that the $\Coker(\beta^{\vee \vee})$ is
  $\Ext_{\mathbb{Z}}^{1}(\dg(\beta), \mathbb{Z})$.  Since
  $\Ext_{\mathbb{Z}}^{1}(\dg(\beta), \mathbb{Z})$ is finite, we see
  that $\beta$ cannot be isomorphic to $\beta^{\vee \vee}$.
  
  Conversely, assume that the cokernel of $\beta$ is finite.  To
  compute the map $\beta^{\vee \vee}$, we first construct a projective
  resolution of $\dg(\beta) = (\mathbb{Z}^{n+r})^{\star} / \image( [ B
  \; Q ]^{\star})$.  Applying the Snake Lemma to the diagram
  \[
  \begin{CD}
    0 @>>> \mathbb{Z}^{r} @>>> \mathbb{Z}^{n+r} @>>> \mathbb{Z}^{n}
    @>>> 0\\
    && @| @VV{[B \; Q]}V @VV{\beta}V \\    
    0 @>>> \mathbb{Z}^{r} @>{Q}>> \mathbb{Z}^{d+r} @>>> N @>>> 0
  \end{CD}
  \]
  establishes that $\Coker([B \; Q]) = \Coker(\beta)$ and $\Ker([B \;
  Q]) = \Ker(\beta)$.  Hence, the complex $0 \rightarrow \Ker(\beta)
  \rightarrow \mathbb{Z}^{n + r} \xrightarrow{[B \; Q]} \mathbb{Z}^{d
    + r} \rightarrow 0$ is a projective resolution of $\Coker(\beta)$.
  Since $\Ext_{\mathbb{Z}}^{i} ( \Coker(\beta),\mathbb{Z} )$ can be
  computed from this resolution and $\Coker(\beta)^{\star} = 0$, we
  deduce that $[B \; Q]^{\star}$ is injective and $0 \rightarrow
  (\mathbb{Z}^{d + r})^{\star} \xrightarrow{[B \; Q]^{\star}}
  (\mathbb{Z}^{n + r})^{\star} \rightarrow 0$ is a projective
  resolution of $\dg(\beta)$.
  
  Since the dual map $\beta^{\vee}$ is the composition of the
  inclusion map $(\mathbb{Z}^{n})^{\star} \rightarrow
  (\mathbb{Z}^{n+r})^{\star}$ and the quotient map
  $(\mathbb{Z}^{n+r})^{\star} \rightarrow \dg(\beta)$, it follows that
  $\dg(\beta^{\vee}) = (\mathbb{Z}^{n + d + r})^{\star \star} / \image
  \left[ 
    \begin{smallmatrix} I_{n} & B^{\star} \\ 0 & Q^{\star}
    \end{smallmatrix} \right]^{\star}$
  and $\beta^{\vee \vee}$ is the composition of inclusion
  $(\mathbb{Z}^{n})^{\star\star} \rightarrow (\mathbb{Z}^{n + d +
    r})^{\star\star}$ and the quotient map $(\mathbb{Z}^{n + d +
    r})^{\star\star} \rightarrow \dg(\beta^{\vee})$.  Because
  $\mathbb{Z}^{m}$ is naturally isomorphic to
  $(\mathbb{Z}^{m})^{\star\star}$, it follows that $\dg(\beta^{\vee})$
  is naturally isomorphic to $(\mathbb{Z}^{d + r} / \image(Q)) = N$
  and $\beta^{\vee \vee}$ is naturally isomorphic to $\beta$.  Lastly,
  our resolution of $\dg(\beta)$ also implies that $H^{0}\big(
  \cone(\beta)^{\star} \big) = 0$ and thus the long exact sequence
  which gives \eqref{eq:dualmap} proves the second part of the
  proposition.
\end{proof}

The operator $(-)^{\vee}$ is also well-behaved in short exact
sequences.

\begin{lemma} \label{lem:ses}
  Given a commutative diagram
  \begin{equation} \label{eq:sesOfMaps}
    \begin{CD}
      0 @>>> \mathbb{Z}^{n_{1}} @>>> \mathbb{Z}^{n_{2}} @>>>
      \mathbb{Z}^{n_{3}} @>>> 0\\
      && @VV{\beta_{1}}V @VV{\beta_{2}}V @VV{\beta_{3}}V\\
      0 @>>> N_{1} @>>> N_{2} @>>> N_{3} @>>> 0
    \end{CD}
  \end{equation}
  in which the rows are exact and the columns have finite cokernels,
  there is a commutative diagram with exact rows:
  \begin{equation} \label{eq:dualses}
    \begin{CD}
      0 @>>> ({\mathbb{Z}}^{n_{3}})^{\star} @>>>
      (\mathbb{Z}^{n_{2}})^{\star} @>>> (\mathbb{Z}^{n_{1}})^{\star}
      @>>> 0\\
      && @VV{\beta_{3}^{\vee}}V @VV{\beta_{2}^{\vee}}V
      @VV{\beta_{1}^{\vee}}V\\
      0 @>>> \dg(\beta_{3}) @>>> \dg(\beta_{2}) @>>> \dg(\beta_{1})
      @>>> 0
    \end{CD} \quad .
  \end{equation}
\end{lemma}

\begin{proof}
  For $1 \leq i \leq 3$, choose $\bm{E}_{i} :=
  (\mathbb{Z}^{n_{i}})[0]$ as a projective resolution of
  $\mathbb{Z}^{n_{i}}$.  Using Lemma~2.2.8 in \cite{Weibel}, the
  bottom row of \eqref{eq:sesOfMaps} lifts to an exact sequence of
  projective resolutions $0 \rightarrow \bm{F}_{1} \rightarrow
  \bm{F}_{2} \rightarrow \bm{F}_{3} \rightarrow 0$.  Hence, the
  diagram \eqref{eq:sesOfMaps} produces a commutative diagram of
  cochain complexes with exact rows:
  \[
  \begin{CD}
    0 @>>> \bm{E}_{1} @>>> \bm{E}_{2} @>>> \bm{E}_{3} @>>> 0\\
    && @VVV @VVV @VVV\\
    0 @>>> \bm{F}_{1} @>>> \bm{F}_{2} @>>> \bm{F}_{3} @>>> 0
  \end{CD} \quad .
  \]
  The naturality of the mapping cone and the functor $(-)^{\star}$
  yield a commutative diagram with exact rows and columns:
  \begin{equation} \label{eq:3*3}
    \begin{CD}
      && 0 && 0 && 0 \\
      && @VVV @VVV @VVV \\
      0 @>>> \bm{E}_{3}[1]^{\star} @>>> \bm{E}_{2}[1]^{\star} @>>>
      \bm{E}_{1}[1]^{\star} @>>> 0 \\
      && @VVV @VVV @VVV \\
      0 @>>> \cone(\beta_{3})^{\star} @>>> \cone(\beta_{2})^{\star}
      @>>> \cone(\beta_{1})^{\star} @>>> 0\\
      && @VVV @VVV @VVV \\      
      0 @>>> \bm{F}_{3}^{\star} @>>> \bm{F}_{2}^{\star} @>>>
      \bm{F}_{1}^{\star} @>>> 0\\
      && @VVV @VVV @VVV \\
      && 0 && 0 && 0 
    \end{CD} \quad .
  \end{equation}
  Since $\Coker(\beta_{i})$ is finite and $\bm{E}_{i} =
  (\mathbb{Z}^{n_{i}})[0]$, both $H^{j}\bigl( \cone(\beta_{i})^{\star}
  \bigr) = 0$ and $H^{j}(\bm{E}_{i}[1]^{\star}) = 0$ for all $j \neq
  1$ and $1 \leq i \leq 3$.  Hence, taking the cohomology of
  \eqref{eq:3*3} yields \eqref{eq:dualses}.
\end{proof}

\section{Toric Deligne-Mumford Stacks} \label{sec:stackytoric}
\setcounter{equation}{0}

The purpose of this section is to associate a smooth Deligne-Mumford
stack to certain combinatorial data.  This construction is inspired by
the quotient construction for toric varieties; for example see
\cite{Cox}.

Let $N$ be a finitely generated abelian group of rank $d$.  We write
$\overline{N}$ for the lattice generated by $N$ in the $d$-dimensional
$\mathbb{Q}$-vector space $N_{\mathbb{Q}} := N \otimes_{\mathbb{Z}}
\mathbb{Q}$.  The natural map $N \rightarrow \overline{N}$ is denoted
by $b \mapsto \bar{b}$.  Let $\Sigma$ be a rational \emph{simplicial}
fan in $N_{\mathbb{Q}}$; every cone $\sigma \in \Sigma$ is generated
by linearly independent vectors.  Let $\rho_{1}, \dotsc, \rho_{n}$ be
the rays (one-dimensional cones) in $\Sigma$.  We assume that
$\rho_{1}, \dotsc, \rho_{n}$ span $N_{\mathbb{Q}}$ and we fix an
element $b_{i} \in N$ such that $\bar b_{i}$ generates the cone
$\rho_{i}$ for $1 \leq i \leq n$.  The set $\{ b_{1}, \dotsc, b_{n}
\}$ defines a homomorphism $\beta \colon \mathbb{Z}^{n} \rightarrow N$
with finite cokernel.  The triple $\bm{\Sigma} := \bigl( N, \Sigma,
\beta \bigr)$ is called a \emph{stacky fan}.

The stacky fan $\bm{\Sigma}$ encodes a group action on a quasi-affine
variety $Z$.  To describe this action, let $\mathbb{C}[z_{1}, \dotsc,
z_{n}]$ be the coordinate ring of $\mathbb{A}^{n}$.  The quasi-affine
variety $Z$ is the open subset defined by the reduced monomial ideal
$J_{\Sigma} := \big\langle \prod_{\rho_{i} \nsubseteq \sigma} z_{i} :
\sigma \in \Sigma \big\rangle$; in other words, $Z := \mathbb{A}^{n} -
\mathbb{V}(J_{\Sigma})$.  The $\mathbb{C}$\nobreakdash-valued points
of $Z$ are the $z \in \mathbb{C}^{n}$ such that the cone generated by
the set $\{ \rho_{i} : \text{$z_{i} = 0$} \}$ belongs to $\Sigma$.  We
equip $Z$ with an action of the group $G := \Hom_{\mathbb{Z}}(
\dg(\beta), \mathbb{C}^{*})$ as follows.  By applying
$\Hom_{\mathbb{Z}}(-, \mathbb{C}^{*})$ to the dual map $\beta^{\vee}
\colon (\mathbb{Z}^{n})^{\star} \rightarrow \dg(\beta)$ (see
Section~\ref{sec:gale}), we obtain a homomorphism $\alpha \colon G
\rightarrow (\mathbb{C}^{*})^{n}$.  The natural action of
$(\mathbb{C}^{*})^{n}$ on $\mathbb{A}^{n}$ induces an action of $G$ on
$\mathbb{A}^{n}$.  Since $\mathbb{V}(J_{\Sigma})$ is a union of
coordinate subspaces, $Z$ is $G$-invariant.

The quotient stack $\mathcal{X}(\bm{\Sigma}) := [Z/G]$ is the Artin
stack associated to the groupoid $s, t \colon Z \times G
\rightrightarrows Z$ where $s$ is the projection onto the first factor
and $t$ is given by the $G$-action on $Z$.  If $S$ is a scheme, then
the objects in $[Z/G](S)$ are principal $G$-bundles $E \rightarrow S$
with a $G$-equivariant map $E \rightarrow Z$ and the morphisms are
isomorphisms which preserve the map to $Z$.  Since $Z$ is smooth,
$\mathcal{X}(\bm{\Sigma})$ is a smooth algebraic stack; see
Remark~10.13.2 in \cite{LaumonMoret}.  The next proposition shows that
$\mathcal{X}(\bm{\Sigma})$ is in fact a Deligne-Mumford stack.  We
call $\mathcal{X}(\bm{\Sigma})$ the \emph{toric Deligne-Mumford stack}
associated to the stacky fan $\bm{\Sigma}$.

\begin{lemma} \label{lem:finiteMap}
  The map $Z \times G \rightarrow Z \times Z$ with $(z,g) \mapsto (z,
  z \cdot g)$ is a finite morphism.
\end{lemma}

\begin{proof}
  The morphism of affine schemes $\alpha \colon G \rightarrow
  (\mathbb{C}^{*})^{n}$ corresponds to the map of rings
  $\mathbb{C}[(\mathbb{Z}^{n})^{\star}] \cong {\mathbb{C}}[t_{1}^{\pm
    1}, \dotsc, t_{n}^{\pm 1}] \rightarrow {\mathbb{C}}[\dg(\beta)]$.
  Since the cokernel of $\beta^{\vee}$ is finite, the ring
  $\mathbb{C}[\dg(\beta)]$ is integral over $\mathbb{C}[t_{1}^{\pm 1},
  \dotsc, t_{n}^{\pm 1}]$ and $G \rightarrow \image(\alpha)$ is a
  finite morphism.  Hence, it suffices to prove that $\xi \colon
  \image(\alpha) \times Z \rightarrow Z \times Z$ is also a finite
  morphism.  Because $\Ker(\beta^{\vee}) \cong N^{\star}$, we have $
  \image(\alpha) = \Spec \bigl( \mathbb{C}[t_{1}^{\pm 1}, \dotsc,
  t_{n}^{\pm 1}] / \langle \prod_{i=1}^{n} t_{i}^{\theta(b_{i})} -1 :
  \theta \in N^{\star} \rangle \bigr)$.
  
  We next show that $\xi \colon \image(\alpha) \times Z \rightarrow Z
  \times Z$ is an affine morphism.  For each $\sigma \in \Sigma$, set
  $z_{\widehat{\sigma}} := \prod_{\rho_{i} \nsubseteq \sigma} z_{i}$
  and let $U_{\sigma} := \mathbb{C}^{n} -
  \mathbb{V}(z_{\widehat{\sigma}})$.  The coordinate ring of the open
  affine subset $U_{\sigma}$ is $\mathbb{C}[z_{1}^{\,}, \dotsc,
  z_{n}^{\,}, z_{\widehat{\sigma}}^{-1}]$ and the collection $\{
  U_{\sigma} : \sigma \in \Sigma \}$ covers $Z$.  Therefore, $\{
  U_{\sigma} \times U_{\sigma'} : \sigma, \sigma' \in \Sigma \}$ is an
  open affine cover of $Z \times Z$ and $U_{\sigma} \times U_{\sigma'}
  = \Spec B_{\sigma, \sigma'}$ where $B_{\sigma, \sigma'} =
  \mathbb{C}[z_{1}^{\,}, \dotsc, z_{n}^{\,},
  z_{\widehat{\sigma}}^{-1}, z_{1}'^{\,}, \dotsc, z_{n}'^{\,},
  (z_{\widehat{\sigma}'}')^{-1}]$.  Since coordinate subspaces are
  $G$\nobreakdash-invariant, $\xi^{-1}(U_{\sigma} \times U_{\sigma'})$
  is the affine set
  \[
  G \times (U_{\sigma} \cap U_{\sigma'}) = \Spec A_{\sigma, \sigma'} =
  \Spec \left( \tfrac{\mathbb{C} \left[ t_{1}^{\pm 1}, \dotsc,
        t_{n}^{\pm 1}, z_{1}^{\,}, \dotsc, z_{n}^{\,},
        z_{\widehat{\sigma}}^{-1}, z_{\widehat{\sigma}'}^{-1}
      \right]}{\left\langle \prod_{i=1}^{n} t_{i}^{\theta(b_{i})} -1
        \; : \; \theta \in N^{\star} \right\rangle} \right) \, .
  \]
  The restriction of $\xi$ to this affine set corresponds to the map
  of rings $\zeta \colon B_{\sigma, \sigma'} \rightarrow A_{\sigma,
    \sigma'}$ given by $z_{i} \mapsto z_{i}$ and $z_{i}' \mapsto
  t_{i}z_{i}$ for $1 \leq i \leq n$.
  
  To prove that $\xi$ is finite, we show that $A_{\sigma, \sigma'}$ is
  a finitely generated $B_{\sigma, \sigma'}$-module.  Clearly, the
  $z_{i} \in A_{\sigma, \sigma'}$ and $(z_{\widehat{\sigma}})^{-1}$
  are integral over $B_{\sigma, \sigma'}$.  Since we have
  \[
  t_{i} = \zeta \Big( (z_{\widehat{\sigma}})^{-1} z_{i}'
  \textstyle\prod_{\substack{\rho_{j} \nsubseteq \sigma' \\ j \neq i}}
  z_{j} \Big) \;\; \text{and} \;\; t_{i}^{-1} = \zeta \Big(
  (z_{\widehat{\sigma}'}')^{-1} z_{i}
  \textstyle\prod_{\substack{\rho_{j} \nsubseteq \sigma' \\ j \neq i}}
  z_{j}' \Big) \, ,
  \]
  both $t_{i}$ for $\bar b_{i} \not\in \sigma$ and $t_{i}^{-1}$ for
  $\bar{b}_{i} \not\in \sigma'$ are integral over $B_{\sigma,
    \sigma'}$.  Thus, $t_{i}^{\pm 1}$ is integral when $\bar b_{i}
  \not\in \sigma \cup \sigma'$.  The Separation Lemma (see Section~1.2
  in \cite{Fulton}) implies there is a $\theta \in N^{\star}$ such
  that $\theta(b_{i}) > 0$ if $\bar b_{i} \in \sigma$ and $\bar b_{i}
  \not\in \sigma'$; $\theta(b_{i}) < 0$ if $\bar b_{i} \not\in \sigma$
  and $\bar b_{i} \in \sigma'$; and $\theta(b_{i}) = 0$ if $\bar b_{i}
  \in \sigma \cap \sigma'$.  Hence, the relation $\prod_{i}
  t_{i}^{\theta(b_{i})} = 1$ can be rewritten as
  $t_{i}^{\theta(b_{i})} = \prod_{j \neq i} t_{j}^{- \theta(b_{j})}$
  and our assumptions on $\theta$ imply that the right hand side is
  integral over $B_{\sigma, \sigma'}$.  It follows that $t_{i}^{\pm
    1}$ is integral over $B_{\sigma, \sigma'}$ when $\bar{b}_{i}
  \not\in \sigma \cap \sigma'$.  Because $\sigma \cap \sigma'$ is
  simplicial, $\bar{b}_{i} \in \sigma \cap \sigma'$ implies that the
  relations $\{ \prod_{i} t_{i}^{\theta(b_{i})} = 1 : \theta \in
  N^{\star} \}$ allow one to express a power of $t_{i}^{\pm 1}$ as a
  product of $t_{j}^{\pm 1}$ for $\bar b_{j} \not\in \sigma \cap
  \sigma'$.  This shows that $t_{i}^{\pm 1}$ for $1 \leq i \leq n$ is
  integral over $B_{\sigma, \sigma'}$.  Lastly, we have
  $(z_{\widehat{\sigma}'})^{-1} = \zeta((z_{\widehat{\sigma}'}')^{-1})
  \prod_{\rho_{i} \nsubseteq \sigma'} t_{i}$ which implies $A_{\sigma,
    \sigma'}$ is integral over $B_{\sigma,\sigma'}$ and completes the
  proof.
\end{proof}

\begin{proposition} 
  The quotient $\mathcal{X}(\bm{\Sigma})$ is a Deligne-Mumford stack.
\end{proposition}

\begin{proof}
  By Corollary~2.2 in \cite{Edidin} (or Example~7.17 in
  \cite{Vistoli}), it is enough to show that the stabilizers of the
  geometric points of $Z$ are finite and reduced.
  Lemma~\ref{lem:finiteMap} shows that the map $Z \times G \rightarrow
  Z \times Z$ defined by $(z, g) \mapsto (z, z \cdot g)$ is a finite
  morphism.  It follows that each stabilizer is a finite group scheme.
  Since we are working in characteristic zero, all finite group
  schemes are reduced.
\end{proof}

\begin{remark}
  In \cite{Lafforgue}, a ``toric stack'' is defined to be the quotient
  of a toric variety by its torus.  Since such a quotient is never a
  Deligne-Mumford stack, $\mathcal{X}(\bm{\Sigma})$ is \emph{not} a
  ``toric stack''.
\end{remark}

\begin{remark}
  The definition of $\mathcal{X}(\bm{\Sigma})$ does not depend on the
  fan $\Sigma$ being simplicial.  However, $\mathcal{X}(\bm{\Sigma})$
  is a Deligne-Mumford stack if and only if the fan $\Sigma$ is
  simplicial.
\end{remark}

As the next example indicates, our construction produces some classic
Deligne-Mumford stacks.

\begin{example}
  Let $\Sigma$ be the complete fan in $\mathbb{Q}$ and consider the
  subset $\bigl\{ (2, 1), (-3, 0) \bigr\}$ of $N := \mathbb{Z} \oplus
  \mathbb{Z}/ 2 \mathbb{Z}$.  This data defines a stacky fan
  $\bm{\Sigma}$.  From Example~\ref{exa:dual}, we know $\beta^{\vee}
  \colon \mathbb{Z}^{2} \rightarrow \dg(\beta) \cong \mathbb{Z}$ is
  given by the matrix $\left[
    \begin{smallmatrix} 6 & 4 \end{smallmatrix} \right]$.
  Furthermore, $Z := \mathbb{A}^{2} - \{ (0,0) \}$ and $\lambda \in G
  \cong \mathbb{C}^{*}$ acts by $(z_{1}, z_{2}) \mapsto (\lambda^{6}
  z_{1}, \lambda^{4} z_{2})$.  In this case,
  $\mathcal{X}(\bm{\Sigma})$ is precisely the moduli stack of elliptic
  curves $\overline{\mathcal{M}}_{1,1}$; see Page~126 in
  \cite{DeligneRapoport}.
\end{example}

To illustrate that a toric Deligne-Mumford stack depends on the set
$\{b_{i}\}$, we include the following:

\begin{example}
  Let $\Sigma$ be the complete fan in $\mathbb{Q}$, which implies $Z
  := \mathbb{A}^{2} - \{ (0,0) \}$, and let $N := \mathbb{Z} \oplus
  \mathbb{Z}/3 \mathbb{Z}$.  If $\beta_{1} \colon \mathbb{Z}^{2}
  \rightarrow N$ corresponds to the set $ \{ (1,0), (-1,1) \}$ and
  $\bm{\Sigma}_{1} = \bigl( N, \Sigma, \beta_{1} \bigr)$, then the
  dual map $\beta_{1}^{\vee} \colon \mathbb{Z}^{2} \rightarrow
  \dg(\beta) \cong \mathbb{Z}$ is given by the matrix $\left[
    \begin{smallmatrix} 3 & 3 \end{smallmatrix} \right]$ 
  and $\lambda \in G_{1} \cong \mathbb{C}^{*}$ acts by $(z_{1}, z_{2})
  \mapsto (\lambda^{3} z_{1}, \lambda^{3} z_{2})$.  On the other hand,
  if $\beta_{2} \colon \mathbb{Z}^{2} \rightarrow N$ corresponds to
  the set $ \{ (1,0), (-1,0) \}$, then $\beta_{2}^{\vee} \colon
  \mathbb{Z}^{2} \rightarrow \dg(\beta) \cong \mathbb{Z} \oplus
  \mathbb{Z}/3 \mathbb{Z}$ is given by $\left[
    \begin{smallmatrix} 1 & 1 \\ 0 & 0 \end{smallmatrix} \right]$ 
  and $(\lambda_{1}, \lambda_{2}) \in G_{2} \cong \mathbb{C}^{*} \times
  \bm{\mu}_{3}$ acts by $(z_{1}, z_{2}) \mapsto (\lambda_{1} z_{1},
  \lambda_{1} z_{2})$.  Therefore, for the stacky fan
  ${\bm{\Sigma}}_{2}=(N,\Sigma,\beta_2)$,
  $\mathcal{X}({\bm{\Sigma}}_{2})$ is the quotient of $\mathbb{P}^{1}$
  by a trivial action of the $\mathbb{Z}/ 3 \mathbb{Z}$ and
  $\mathcal{X}(\bm{\Sigma}_{1}) \not\cong
  \mathcal{X}(\bm{\Sigma}_{2})$.
\end{example}

The last result in this section makes the relationship between toric
Deligne-Mumford stacks and toric varieties more explicit.  Recall that
a \emph{coarse moduli space} of a Deligne-Mumford stack $\mathcal{X}$
is an algebraic space $X$ with a morphism $\pi \colon \mathcal{X}
\rightarrow X$ such that
\begin{itemize}
\item for all algebraically closed fields $k$, the map $\pi(k) \colon
  \mathcal{X}(k) \rightarrow X(k)$ is a bijection;
\item given any algebraic space $X'$ and any morphism $\pi' \colon
  \mathcal{X} \rightarrow X'$, there is a unique morphism $\chi \colon
  X \rightarrow X'$ such that $\pi' = \chi \circ \pi$.
\end{itemize}

\begin{proposition} \label{pro:coarse}
  The toric variety $X(\Sigma)$ is the coarse moduli space of
  $\mathcal{X}(\bm{\Sigma})$.
\end{proposition}

\begin{proof}
  By Proposition~4.2 in \cite{Edidin}, it is enough to show that the
  toric variety $X(\Sigma)$ is the universal geometric quotient of $Z$
  by $G$.  Under the additional assumptions that $N = \overline{N}$
  and that the $b_{i} = \bar b_{i}$ are the unique minimal lattice
  points generating the rays in $\Sigma$, this is Theorem~2.1 in
  \cite{Cox}.  The reader can verify that the proof presented in
  \cite{Cox} extends to our situation without any significant changes.
\end{proof}

\section{Closed and Open Substacks} \label{sec:openclosed}
\setcounter{equation}{0}

This section explains how the stacky fan $\bm{\Sigma}$ encodes certain
closed and open substacks of $\mathcal{X}(\bm{\Sigma})$.  We also
express the inertia stack $\mathcal{I} \bigl( \mathcal{X}(\bm{\Sigma})
\bigr)$ as a disjoint union of certain closed substacks.

To describe the connection between the combinatorics of the stacky fan
$\bm{\Sigma}$ and the substacks of $\mathcal{X}(\bm{\Sigma})$, we use
the theory of groupoids; see \cite{Moerdijk} for an introduction.
Recall that a homomorphism of groupoids $\Theta \colon (R'
\rightrightarrows U') \longrightarrow (R \rightrightarrows U)$ is
called a \emph{Morita equivalence} if
\begin{enumerate}
\item the square
\[
\begin{CD}
  R' @>{(s,t)}>> U' \times U' \\
  @V{\Theta}VV @VV{\Theta \times \Theta}V\\
  R @>{(s,t)}>> U \times U
\end{CD}
\]
is Cartesian, and
\item the morphism $t \circ \pr_{1} \colon U' \times_{\Theta, U, s} R
  \rightarrow U$ is locally surjective.  In other words, $U$ has an
  open covering $\{ U_{i} \rightarrow U \}$ in the \'{e}tale topology
  such that each $U_{i} \rightarrow U$ factors through $U$.
\end{enumerate}
The key observation is that two groupoids are Morita equivalent if and
only if the associated stacks are isomorphic.

Fix a cone $\sigma$ in the fan $\Sigma$.  Let $N_{\sigma}$ be the
subgroup of $N$ generated by the set $\{ b_{i} : \rho_{i} \subseteq
\sigma \}$ and let $N(\sigma)$ be the quotient group $N / N_{\sigma}$.
By extending scalars, the quotient map $N \rightarrow N(\sigma)$
becomes the surjection $N_{\mathbb{Q}} \rightarrow
N(\sigma)_{\mathbb{Q}}$.  The quotient fan $\Sigma / \sigma$ in
$N(\sigma)_{\mathbb{Q}}$ is the set $\{ \tilde{\tau} = \tau +
(N_{\sigma})_{\mathbb{Q}} : \text{$\sigma \subseteq \tau$ and $\tau
  \in \Sigma$} \}$ and the link of $\sigma$ is the set $\link(\sigma)
:= \{ \tau : \tau + \sigma \in \Sigma, \tau \cap \sigma = 0 \}$.  For
each ray $\rho_{i}$ in $\link(\sigma)$, we write $\tilde{\rho}_{i}$
for the ray in $\Sigma / \sigma$ and $\tilde{b}_{i}$ for the image of
$b_{i}$ in $N(\sigma)$.  To ensure that the quotient fan satisfies our
hypothesis for constructing toric Deligne-Mumford stacks, we require
the following:

\begin{condition} \label{con:spanQuot}
  The rays $\tilde{\rho}_{i}$ span $N(\sigma)_{\mathbb{Q}}$.
\end{condition}

\noindent
Note that if $\Sigma$ is a complete fan, then every cone $\sigma$
satisfies Condition~\ref{con:spanQuot}.

Let $\ell$ be the number of rays in $\link(\sigma)$ and let
$\beta(\sigma) \colon \mathbb{Z}^{\ell} \rightarrow N(\sigma)$ be the
map determined by the set $\{ \tilde{b}_{i} : \rho_{i} \in
\link(\sigma) \}$.  The quotient stacky fan $\bm{\Sigma / \sigma}$ is
the triple $\bigl( N(\sigma), \Sigma / \sigma, \beta(\sigma) \bigr)$.

\begin{proposition} \label{pro:closed}
  If $\sigma$ is a cone in the stacky fan $\bm{\Sigma}$ which
  satisfies Condition~$\ref{con:spanQuot}$, then
  $\mathcal{X}(\bm{\Sigma / \sigma})$ defines a closed substack of
  $\mathcal{X}(\bm{\Sigma})$.
\end{proposition}

\begin{proof}
  By definition, $\mathcal{X}(\bm{\Sigma})$ is $[Z/G]$.  Let
  $W(\sigma)$ be the closed subvariety of $Z$ defined by the ideal
  $J(\sigma) := \langle z_{i}: \rho_{i} \subseteq \sigma \rangle$ in
  $\mathbb{C}[z_{1}, \dotsc, z_{n}]$.  The $\mathbb{C}$-valued points
  of $W(\sigma)$ are the $z \in \mathbb{C}^{n}$ such that the cone
  spanned by $\{ \rho_{i} : z_{i} = 0 \}$ contains $\sigma$ and
  belongs to $\Sigma$.  Hence, $\rho_{i} \not\subseteq \sigma \cup
  \link(\sigma)$ implies that $z_{i} \neq 0$.  Since $J(\sigma)$
  defines a coordinate subspace, $W(\sigma)$ is
  $G$\nobreakdash-invariant and the groupoid $W(\sigma) \times G
  \rightrightarrows W(\sigma)$ defines a closed substack of
  $\mathcal{X}(\bm{\Sigma})$.  It remains to show that
  $\mathcal{X}(\bm{\Sigma / \sigma})$ is the stack associated to
  $W(\sigma) \times G \rightrightarrows W(\sigma)$.
  
  To begin, we construction a homomorphism from $W(\sigma) \times G
  \rightrightarrows W(\sigma)$ to the defining groupoid of
  $\mathcal{X}(\bm{\Sigma / \sigma})$.  By renumbering the $\rho_{i}$,
  we may assume that $\tilde{\rho}_{1}, \dotsc, \tilde{\rho}_{\ell}$
  are the rays in $\link(\sigma)$.  If $\mathbb{C}[\tilde{z}_{1},
  \dotsc, \tilde{z}_{\ell}]$ is the coordinate ring of
  $\mathbb{A}^{\ell}$, then 
  \[
  J_{\Sigma / \sigma} := \bigl\langle \textstyle\prod_{\rho_{i}
    \nsubseteq \tau} \tilde{z}_{i} : \text{$\sigma \subseteq \tau$ and
    $\tau \in \Sigma$} \bigr\rangle \, .
  \]  
  By definition, $\mathcal{X}(\bm{\Sigma / \sigma}) :=
  [Z(\sigma)/G(\sigma)]$ where $Z(\sigma) := \mathbb{A}^{\ell} -
  \mathbb{V}(J_{\Sigma / \sigma})$ and $G(\sigma)$ is the group
  $\Hom_{\mathbb{Z}}(\dg(\beta(\sigma)), \mathbb{C}^{*})$.  Let $m :=
  \dim \sigma$.  The description of the $\mathbb{C}$-valued points of
  $W(\sigma)$ shows the projection $\mathbb{A}^{n} \rightarrow
  {\mathbb{A}}^{\ell}$ induces a surjection $\varphi_{0} \colon
  W(\sigma) \rightarrow Z(\sigma)$ with $\Ker(\varphi_{0}) =
  (\mathbb{C}^{*})^{n-\ell - m}$.  Applying Lemma~\ref{lem:ses} to the
  commutative diagram
  \[
  \begin{CD}
    0 @>>> \mathbb{Z}^{n - \ell} @>>> \mathbb{Z}^{n} @>>>
    \mathbb{Z}^{\ell} @>>> 0\\
    && @VV{\tilde{\beta}}V @VV{\beta}V @VV{\beta(\sigma)}V \\
    0 @>>> N_{\sigma} @>>> N @>>> N(\sigma) @>>> 0
  \end{CD}
  \]
  produces the commutative diagram with exact rows
  \begin{equation} \label{eq:closedDual}
    \begin{CD}
      0 @>>> (\mathbb{Z}^{\ell})^{\star} @>>> (\mathbb{Z}^{n})^{\star}
      @>>> (\mathbb{Z}^{n-\ell})^{\star} @>>> 0\\
      && @VV{\beta(\sigma)^{\vee}}V @VV{\beta^{\vee}}V
      @VV{\tilde{\beta}^{\vee}}V\\
      0 @>>> \dg(\beta(\sigma)) @>>> \dg(\beta) @>>>
      \dg(\tilde{\beta}) @>>> 0
    \end{CD} \quad .
  \end{equation}
  Since the cone $\sigma$ is simplicial, $N_{\sigma} \cong
  \mathbb{Z}^{m}$ and $\dg(\tilde{\beta}) \cong \mathbb{Z}^{n-\ell -
    m}$.  Applying the functor $\Hom_{\mathbb{Z}}(-, \mathbb{C}^{*})$
  to \eqref{eq:closedDual} gives the diagram with split exact rows
  \[ 
  \begin{CD}
    0 @>>> (\mathbb{C}^{*})^{n-\ell - m} @>>> G @>{\varphi_{1}}>>
    G(\sigma) @>>> 0\\
    && @VVV @VV{\alpha}V @VV{\alpha(\sigma)}V\\
    0 @>>> (\mathbb{C}^{*})^{n-\ell} @>>> (\mathbb{C}^{*})^{n} @>>>
    (\mathbb{C}^{*})^{\ell} @>>> 0
  \end{CD} \quad .
  \] 
  Hence, $\Phi := (\varphi_{0} \times \varphi_{1}, \varphi_{0})$ is a
  homomorphism of groupoids from $W(\sigma) \times G \rightrightarrows
  W(\sigma)$ to $Z(\sigma) \times G(\sigma) \rightrightarrows
  Z(\sigma)$.
  
  To prove that $\mathcal{X}(\bm{\Sigma / \sigma})$ is the stack
  associated to $W(\sigma) \times G \rightrightarrows W(\sigma)$, it
  suffices to show that $\Phi$ is a Morita equivalence.  First, the
  commutative diagram
  \[
  \begin{CD}
    Z(\sigma) \times G(\sigma) \times (\mathbb{C}^{*})^{2(n-\ell-m)}
    @<{\cong}<< W(\sigma) \times G @>{\varphi_{0} \times
      \varphi_{1}}>> Z(\sigma) \times G(\sigma)\\
    @VV{(s,t,\text{id})}V @VV{(s,t)}V @VV{(s,t)}V\\
    Z(\sigma) \times Z(\sigma) \times (\mathbb{C}^{*})^{2(n-\ell-m)}
    @<{\cong}<< W(\sigma) \times W(\sigma) @>{\varphi_{0} \times
      \varphi_{0}}>> Z(\sigma) \times Z(\sigma)
  \end{CD}
  \]
  shows that $W(\sigma) \times G = \bigl( Z(\sigma) \times G(\sigma)
  \bigr) \times_{\varphi_{0} \times \varphi_{0}, Z(\sigma) \times
    Z(\sigma), (s,t)} \bigl( W(\sigma) \times W(\sigma) \bigr)$.
  Second, we have $\bigl( Z(\sigma) \times G(\sigma) \bigr) \times_{s,
    Z(\sigma), \varphi_{0}} W(\sigma) \cong Z(\sigma) \times G(\sigma)
  \times \mathbb{C}^{n-\ell-m}$ which implies that the map $t \circ
  \pi_{1} \colon \bigl( Z(\sigma) \times G(\sigma) \bigr) \times_{s,
    Z(\sigma), \varphi_{0}} W(\sigma) \rightarrow Z(\sigma)$ splits.
  Therefore, $\Phi$ is a Morita equivalence and
  $\mathcal{X}(\bm{\Sigma / \sigma})$ defines a closed substack of
  $\mathcal{X}(\bm{\Sigma})$.
\end{proof}

Viewing $\sigma \in \Sigma$ as the fan consisting of the cone $\sigma$
and all its faces, we can identify $\sigma$ with an open substack of
$\mathcal{X}(\bm{\Sigma})$.  This substack has a particularly nice
description when $\sigma$ is of maximal dimension; $\dim \sigma = d :=
\rank N$.  In this case, let $\beta_{\sigma} \colon \mathbb{Z}^{d}
\rightarrow N$ be the map determined by the set $\{ b_{i} : \rho_{i}
\subseteq \sigma \}$.  The induced stacky fan $\bm{\sigma}$ is the
triple $\bigl( N, \sigma, \beta_{\sigma} \bigr)$.

\begin{proposition} \label{pro:open}
  If $\sigma$ is a $d$-dimensional cone in the stacky fan
  $\bm{\Sigma}$, then $\mathcal{X}(\bm{\sigma})$ defines an open
  substack of $\mathcal{X}(\bm{\Sigma})$.  Moreover,
  $\mathcal{X}(\bm{\sigma})$ is isomorphic the quotient of
  $\mathbb{C}^{d}$ by the finite abelian group $N(\sigma)$.
\end{proposition}

\begin{proof}
  As in Lemma~\ref{lem:finiteMap}, let $U_{\sigma}$ be the open
  subvariety of $Z$ defined by the monomial $z_{\widehat{\sigma}} :=
  \prod_{\rho_{i} \nsubseteq \sigma} z_{i}$.  The $\mathbb{C}$-valued
  points of $U_{\sigma}$ are the $z \in \mathbb{C}^{n}$ such that for
  each $z_{i} = 0$ the ray $\rho_{i}$ is contained in $\sigma$.  Since
  $\mathbb{V}(z_{\widehat{\sigma}})$ is a union of coordinate
  subspaces, $U_{\sigma}$ is $G$-invariant and the groupoid
  $U_{\sigma} \times G \rightrightarrows U_{\sigma}$ defines an open
  substack of $\mathcal{X}(\bm{\Sigma})$.  It remains to show that
  $\mathcal{X}(\bm{\sigma})$ is the stack associated to $U_{\sigma}
  \times G \rightrightarrows U_{\sigma}$.
  
  We construct a homomorphism from the defining groupoid of
  $\mathcal{X}(\bm{\sigma})$ to $U_{\sigma} \times G \rightrightarrows
  U_{\sigma}$.  Since $\sigma$ is a $d$-dimensional simplicial cone,
  $J_{\sigma} = \langle 1 \rangle$ and $Z_{\sigma} := \mathbb{A}^{d}$.
  By definition, $\mathcal{X}(\bm{\sigma}) := [Z_{\sigma} /
  G_{\sigma}]$ where $G_{\sigma} := \Hom_{\mathbb{Z}} (
  \dg(\beta_{\sigma}), \mathbb{C}^{*})$.  The description of the
  $\mathbb{C}$-valued points of $U_{\sigma}$ yields a closed embedding
  $\psi_{0} \colon Z_{\sigma} \rightarrow U_{\sigma}$ where
  \[
  \psi_{0}(Z_{\sigma}) = \mathbb{C}^{d} \times \bm{1} \subset
  \mathbb{C}^{d} \times (\mathbb{C}^{*})^{n-d} \cong U_{\sigma} \, .
  \]
  Applying Lemma~\ref{lem:ses} and the functor $\Hom_{\mathbb{Z}}(-,
  \mathbb{C}^{*})$ to
  \[
  \begin{CD}
    0 @>>> \mathbb{Z}^{d} @>>> \mathbb{Z}^{n} @>>> \mathbb{Z}^{n-d}
    @>>> 0\\
    && @VV{\beta_{\sigma}}V @VV{\beta}V @VVV\\
    0 @>>> N @>>{\text{id}}> N @>>> 0 
  \end{CD}
  \]
  produces the commutative diagram:
  \begin{equation} \label{eq:openDual}
    \begin{CD}
      0 @>>> G_{\sigma} @>>> G @>{\psi_{1}}>> (\mathbb{C}^{*})^{n-d}
      @>>> 0\\
      && @VV{\alpha_{\sigma}}V @VV{\alpha}V @VV{\text{id}}V\\
      0 @>>> (\mathbb{C}^{*})^{d} @>>> (\mathbb{C}^{*})^{n} @>>>
      (\mathbb{C}^{*})^{n-d} @>>> 0
    \end{CD} \quad .
  \end{equation}
  Hence, $\Psi := (\psi_{0} \times \psi_{1}, \psi_{0})$ is a
  homomorphism of groupoids from $Z_{\sigma} \times G_{\sigma}
  \rightrightarrows Z_{\sigma}$ to $U_{\sigma} \times G
  \rightrightarrows U_{\sigma}$ and an element $g \in G$ belongs to
  $G_{\sigma}$ if and only if $(Z_{\sigma} \cdot g) \cap Z_{\sigma}
  \neq \emptyset$.
  
  Next, we establish that $G_{\sigma} \cong N(\sigma)$.  The
  definition of $N(\sigma)$ gives the exact sequence
  \[
  0 \longrightarrow \mathbb{Z}^{d+r} \xrightarrow{\;\; [B_{\sigma} \;
    Q]\;\; } \mathbb{Z}^{d+r} \longrightarrow N(\sigma)
  \longrightarrow 0
  \] 
  where $B_{\sigma}$ is the submatrix of $B$ whose columns correspond
  to the $\rho_{i} \subseteq \sigma$.  Since $N(\sigma)^{\star} = 0$,
  we obtain the exact sequence
  \[
  0 \longrightarrow (\mathbb{Z}^{d+r})^{\star} \xrightarrow{\;\;
    [B_{\sigma} \; Q]^{\star}} (\mathbb{Z}^{d+r})^{\star}
  \longrightarrow \Ext_{\mathbb{Z}}^{1} \bigl( N(\sigma), \mathbb{Z}
  \bigr) \longrightarrow 0
  \]
  which implies that $\dg(\beta) = \Ext_{\mathbb{Z}}^{1} \bigl(
  N(\sigma), \mathbb{Z} \bigr) = \Hom_{\mathbb{Z}} \bigl( N(\sigma),
  \mathbb{Q}/\mathbb{Z} \bigr)$.  Hence, the group $G_{\sigma}$ is
  $\Hom_{\mathbb{Z}} \bigl( \Hom_{\mathbb{Z}} \bigl( N(\sigma),
  \mathbb{Q} / \mathbb{Z} \bigr), \mathbb{C}^{*} \bigr)$.  We identify
  $\mathbb{Q}/ \mathbb{Z}$ with a subgroup of $\mathbb{C}^{*}$ via the
  map $p \mapsto \exp(2 \pi \sqrt{-1} p)$ to obtain a natural
  homomorphism from $N(\sigma)$ to $G_{\sigma}$.  By expressing
  $N(\sigma)$ as a direct sum of cyclic groups, one verifies that this
  map is an isomorphism.
  
  Finally, to prove that $\mathcal{X}(\bm{\sigma})$ is the stack
  associated to $U_{\sigma} \times G \rightrightarrows U_{\sigma}$, it
  suffices to show that $\Psi$ is a Morita equivalence.  First,
  because an element $g \in G$ belongs to $G_{\sigma}$ if and only if
  $(Z_{\sigma} \cdot g) \cap Z_{\sigma} \neq \emptyset$, the
  commutative diagram
  \[
  \begin{CD}
    Z_{\sigma} \times G_{\sigma} @>{\psi_{0} \times \psi_{1}}>>
    U_{\sigma} \times G\\
    @VV{(s,t)}V @VV{(s,t)}V\\
    Z_{\sigma} \times Z_{\sigma} @>{\psi_{0} \times \psi_{0}}>>
    U_{\sigma} \times U_{\sigma}
  \end{CD}
  \]
  establishes that $Z_{\sigma} \times G_{\sigma} = \bigl( Z_{\sigma}
  \times Z_{\sigma} \bigr) \times_{\psi_{0} \times \psi_{0},
    Z_{\sigma} \times Z_{\sigma}, (s,t)} \bigl( U_{\sigma} \times G
  \bigr)$.  Secondly, we have $\bigl( U_{\sigma} \times G \bigr)
  \times_{s, U_{\sigma}, \psi_{0}} Z_{\sigma} \cong Z_{\sigma} \times
  G$ which implies that $\pi_{1} \colon \bigl( U_{\sigma} \times G
  \bigr) \times_{s, U_{\sigma}, \psi_{0}} Z_{\sigma} \rightarrow
  U_{\sigma} \times G$ corresponds to the closed immersion $\psi_{0}
  \times \text{id} \colon Z_{\sigma} \times G \rightarrow U_{\sigma}
  \times G$.  Lemma~\ref{lem:finiteMap} implies that $t \colon
  U_{\sigma} \times G \rightarrow U_{\sigma}$ is finite.  Since the
  action of $\Coker(\psi_{1})$ on $\psi_{0}(Z_{\sigma})$ surjects onto
  $U_{\sigma}$, we deduce that $t \circ \pi_{1} \colon \bigl(
  U_{\sigma} \times G \bigr) \times_{s, U_{\sigma}, \psi_{0}}
  Z_{\sigma} \rightarrow U_{\sigma}$ is a finite surjective morphism
  of nonsingular varieties and hence flat.  Because the geometric
  fibers of $t \circ \pi_{1}$ correspond to $G_{\sigma}$, a finite set
  of reduced points, the map $t \circ \pi_{1}$ is also \'{e}tale and
  therefore locally surjective.  We conclude that $\Psi$ is a Morita
  equivalence and $\mathcal{X}(\bm{\sigma})$ defines an open substack
  of $\mathcal{X}(\bm{\Sigma})$.
\end{proof}

\begin{remark} \label{rem:covered}
  Assuming that every cone in $\Sigma$ is contained in a
  $d$-dimensional cone, Proposition~\ref{pro:open} produces an
  \'{e}tale atlas of $\mathcal{X}(\bm{\Sigma})$.
\end{remark}

\begin{remark}
  More generally, if $\bm{\Sigma}' := (N',\Sigma', \beta')$ and
  $\bm{\Sigma} := (N, \Sigma, \beta)$ are two stacky fans, then a
  \emph{morphism of stacky fans} is a homomorphism $\phi \colon N'
  \rightarrow N$ satisfying:
  \begin{itemize}
  \item for each cone $\sigma' \in \Sigma'$, there exists a $\sigma
    \in \Sigma$ such that $\phi_{\mathbb{Q}}(\sigma') \subseteq
    \sigma$ where $\phi_{\mathbb{Q}} \colon N' \otimes_{\mathbb{Z}}
    \mathbb{Q} \rightarrow N \otimes_{\mathbb{Z}} \mathbb{Q}$;
  \item the element $\phi(b_{i}')$ is an integer combination of the
    $b_{j} \in N$ where $\bar{b}_{i}' \in \sigma'$, $\bar{b}_{j} \in
    \sigma$ and $\sigma \in \Sigma$ is any cone that contains
    $\phi_{\mathbb{Q}}(\sigma')$.
  \end{itemize}
  For each morphism $\phi \colon \bm{\Sigma}' \rightarrow
  \bm{\Sigma}$, there is a morphism $\mathcal{X}(\bm{\Sigma}')
  \rightarrow \mathcal{X}(\bm{\Sigma})$.  Since we do not make use of
  this construction, the proof is left to the reader.
\end{remark}

For each $d$-dimensional cone $\sigma$ in the stacky fan
$\bm{\Sigma}$, we define $\BBox(\sigma)$ to be the set of elements $v
\in N$ such that $\bar{v} = \sum_{\rho_{i} \subseteq \sigma} q_{i}
\bar{b}_{i}$ for some $0 \leq q_{i} < 1$.  Hence, the set
$\BBox(\sigma)$ is in one-to-one correspondence with the elements in
the finite group $N(\sigma)$.  Let $\BBox(\bm{\Sigma})$ be the union
of $\BBox(\sigma)$ for all $d$-dimensional cones $\sigma \in \Sigma$.
For each $v \in N$, we write $\sigma(\bar{v})$ for the unique minimal
cone containing $\bar{v}$.

\begin{lemma} \label{lem:boxFixed}
  If $\Sigma$ is a complete fan, then the elements $v \in
  \BBox(\bm{\Sigma})$ are in one-to-one correspondence with elements
  $g \in G$ which fix a point of $Z$.  Moreover, we have $[Z^{g}/G]
  \cong \mathcal{X} \bigl( \bm{\Sigma/ \sigma}(\bar{v}) \bigr)$.
\end{lemma}

\begin{proof}
  By definition, an element $v \in \BBox(\bm{\Sigma})$ corresponds to
  an element in $N(\tau)$ for some $d$-dimensional cone $\tau \in
  \Sigma$.  In the proof of Proposition~\ref{pro:open}, we give an
  isomorphism between $N(\tau)$ and and $G_{\tau}$.  Hence, there is a
  bijection sending $v$ to an element $g$ in the subgroup $G_{\tau}
  \subseteq G$.  In addition, \eqref{eq:openDual} implies that $g$ act
  trivially on points $z \in Z$ with $z_{i} = 0$ for all $\rho_{i}
  \subseteq \tau$ which shows that $g$ fixes a point in $Z$.
  
  Conversely, suppose $g \in G$ fixes a point $z \in Z$.  Since the
  action of $G$ on $Z$ is defined via the map $\alpha \colon G
  \rightarrow (\mathbb{C}^{*})^{n}$ where $g \mapsto \bigl(
  \alpha_{1}(g), \dotsc, \alpha_{n}(g) \bigr)$, we see that either
  $\alpha_{i}(g) = 1$ or $z_{i} = 0$ for all $1 \leq i \leq n$.  The
  definition of $Z$ guarantees that there exists a cone in $\Sigma$
  containing all the rays $\rho_{i}$ for which $z_{i} = 0$.  Let
  $\sigma$ be the minimal cone with this property.  Because $\Sigma$
  is a simplicial fan, the ray $\rho_{i}$ is contained in $\sigma$ if
  and only if $\alpha_{i}(g) \neq 1$.  Thus, the closed subvariety
  $W(\sigma)$ defined in Proposition~\ref{pro:closed} is equal to the
  invariant subvariety $Z^{g}$.  Moreover, our choice of $\sigma$
  implies that the element $g$ stabilizes $\psi_{0}(Z_{\tau})$ for
  every $d$-dimensional cone $\tau$ which contains $\sigma$.  It
  follows that $g$ corresponds to an element $v \in
  \BBox(\bm{\Sigma})$.  Finally, $\sigma$ is clearly the intersection
  of all maximal cones $\tau$ for which $v$ corresponds to an element
  in $N(\tau)$.  Therefore, $\sigma = \sigma(\bar{v})$ and
  Proposition~\ref{pro:closed} establishes that $[Z^g/G] =
  [W(\sigma)/G] \cong \mathcal{X} \bigl( \bm{\Sigma/ \sigma}(\bar v)
  \bigr)$.
\end{proof}

For a Deligne-Mumford stack $\mathcal{X}$, its \emph{inertia stack}
$\mathcal{I}(\mathcal{X})$ is defined to be the fibered product
$\mathcal{X} \times_{\Delta, \mathcal{X} \times \mathcal{X}, \Delta}
\mathcal{X}$ where $\Delta$ denotes the diagonal map.  For a scheme
$S$, an object in $\mathcal{I} \bigl( \mathcal{X} \bigr) (S)$ can be
identified with pair $(x, \phi)$ where $x$ is an object in
$\mathcal{X}(S)$ and $\phi$ is an automorphism of $x$.  A morphism
from $(x, \phi) \rightarrow (x', \phi')$ is a morphism $\gamma \colon
x \rightarrow x'$ in $\mathcal{X}(S)$ such that $\gamma \circ \phi =
\phi' \circ \gamma$.  Since we are working over $\mathbb{C}$, the
inertia stack $\mathcal{I}(\mathcal{X})$ is naturally isomorphic to
the stack of representable morphisms from constant cyclotomic gerbes
to $\mathcal{X}$; see Section~4.4 in \cite{AGV}.

\begin{proposition} \label{pro:inertia}
  If $\Sigma$ is a complete fan, then $\mathcal{I} \bigl(
  \mathcal{X}(\bm{\Sigma}) \bigr) = \coprod_{v \in \BBox(\bm{\Sigma})}
  \mathcal{X} \bigl( \bm{\Sigma / \sigma}(\bar{v}) \bigr)$ where
  $\sigma(\bar{v})$ is the minimal cone in $\Sigma$ containing
  $\bar{v}$.
\end{proposition}

\begin{proof}
  Let $S$ be a connected scheme.  An object $x$ of
  $\mathcal{X}(\bm{\Sigma})(S)$ is a principal $G$-bundle $E
  \rightarrow S$ with a $G$-equivariant morphism $f \colon E
  \rightarrow Z$.  An automorphism $\phi$ is an automorphism of the
  principal $G$-bundle $E \rightarrow S$ that is compatible with $E
  \rightarrow Z$.  Since $S$ is connected, $\phi$ corresponds to
  multiplication by an element $g \in G$.  Moreover, because $f$ is
  $G$-equivariant and $f = f \circ \phi$, the map $f$ factors through
  $Z^{g}$.  Hence, the principal $G$-bundle $E \rightarrow S$ with $E
  \rightarrow Z^{g}$ is an object in $[Z^{g} / G](S)$.
  
  For an arbitrary scheme $S$ and an object in $\mathcal{I} \bigl(
  \mathcal{X}(\bm{\Sigma}) \bigr)(S)$, we can assign an object in
  $\coprod_{g \in G} [Z^{g}/G](S)$ by considering the connected
  components of $S$.  Finally, Lemma~\ref{lem:boxFixed} shows that
  $Z^{g} \neq \emptyset$ if and only if $g$ corresponds to an element
  $v \in \BBox(\bm{\Sigma})$ and that $[Z^g/G] \cong \mathcal{X}
  \bigl( \bm{\Sigma/ \sigma}(\bar v) \bigr)$.
\end{proof}

\begin{remark}
  By combining Proposition~\ref{pro:coarse} and
  Proposition~\ref{pro:inertia}, we see that the coarse moduli space
  of $\mathcal{I} \bigl( \mathcal{X}(\bm{\Sigma}) \bigr)$ is
  isomorphic to the disjoint union of $X \big( \Sigma / \sigma(\bar
  v)\big)$ for all $v \in \BBox(\Sigma)$. In particular, we recover
  the description of the twisted sectors in Section~6 of
  \cite{Poddar}.
\end{remark}

\section{Module Structure on $A_{orb}^{*} \bigl(
  \mathcal{X}(\bm{\Sigma}) \bigr)$} \label{sec:chowanddef}
\setcounter{equation}{0}

The goal of this section to describe the orbifold Chow ring of a
complete toric Deligne-Mumford stack as an abelian group.  Throughout
this section, we assume all fans are complete and simplicial and all
Chow rings have rational coefficients.

We first introduce the \emph{deformed group ring}
$\mathbb{Q}[N]^{\bm{\Sigma}}$ associated to the stacky fan
$\bm{\Sigma} = \big(N, \Sigma, \beta \big)$.  As a vector space,
$\mathbb{Q}[N]^{\bm{\Sigma}}$ is simply the group ring
$\mathbb{Q}[N]$; in other words, $ \mathbb{Q}[N]^{\bm{\Sigma}} =
\bigoplus_{c \in N} \mathbb{Q} \cdot y^{c} $ where $y$ is a formal
variable.  Multiplication in $\mathbb{Q}[N]^{\bm{\Sigma}}$ is defined
as follows:
\begin{equation} \label{eq:definingProd}
y^{c_{1}} \cdot y^{c_{2}} :=
\begin{cases}
  y^{c_{1}+c_{2}} & \text{if there exists $\sigma \in \Sigma$ such
    that $\bar c_{1} \in \sigma$ and $\bar c_{2} \in
    \sigma$;} \\
  0 & \text{otherwise.}
\end{cases}
\end{equation}
We endow $\mathbb{Q}[N]^{\bm{\Sigma}}$ with a $\mathbb{Q}$-grading as
follows: if $\bar{c} = \sum_{\rho_{i} \subseteq \sigma(\bar{c})} m_{i}
\bar{b}_{i}$ where $\sigma(\bar{c})$ is the minimal cone in $\Sigma$
containing $\bar{c}$, then $\deg(y^{c}) := \sum m_{i} \in \mathbb{Q}$.

Given a stacky fan $\bm{\Sigma}$, we denote by $S_{\bm{\Sigma}}$ the
subring of $\mathbb{Q}[N]^{\bm{\Sigma}}$ generated over $\mathbb{Q}$
by the monomials $y^{b_i}$.  Since $\Sigma$ is simplicial, the ring
$S_{\bm{\Sigma}}$ is isomorphic to the quotient $\mathbb{Q}[x_{1},
\dotsc, x_{n}]/I_{\Sigma}$ where the ideal $I_{\Sigma}$ is generated
by the square-free monomials $x_{i_{1}}x_{i_{2}} \dotsb x_{i_{s}}$ with
$\rho_{i_{1}} + \dotsb + \rho_{i_{s}} \not\in \Sigma$.  In particular,
$S_{\bm{\Sigma}}$ is a $\mathbb{Z}$-graded ring and $I_{\Sigma}$ is
the Stanley-Reisner ideal associated to $\Sigma$.

To describe the Chow ring of $\mathcal{X}(\bm{\Sigma})$, we need
certain line bundles corresponding to the rays $\rho_{1}, \dotsc,
\rho_{n}$.  Since the category of coherent sheaves on
$\mathcal{X}(\bm{\Sigma})$ is equivalent to the category of
$G$-equivariant sheaves on $Z$ (Example~7.21 in \cite{Vistoli}), we
can define $L_{i}$ for $1 \leq i \leq n$ to be the line bundle on
$\mathcal{X}(\bm{\Sigma})$ corresponding to the trivial line bundle
$\mathbb{C} \times Z$ on $Z$ with the $G$-action on $\mathbb{C}$ is
given by the $i$-th component $\alpha_{i}$ of $\alpha \colon G
\rightarrow (\mathbb{C}^{*})^{n}$.

We first calculate the non-orbifold Chow ring of
$\mathcal{X}(\bm{\Sigma})$.

\begin{lemma} \label{lem:basicChow}
  If $\mathcal{X}(\bm{\Sigma})$ is a complete toric Deligne-Mumford
  stack, then there is an isomorphism of $\mathbb{Z}$-graded rings
  \[
  \frac{S_{\bm{\Sigma}}} {\left\langle \sum_{i=1}^{n} \theta( \bar
      b_{i} ) \cdot y^{b_{i}} : \theta \in N^{\star} \right\rangle}
  \longrightarrow A^{*}\big( \mathcal{X}(\bm{\Sigma}) \big)
  \]
  defined by $y^{b_{i}} \mapsto c_{1}(L_{i})$.
\end{lemma}

\begin{proof}
  For $1 \leq i \leq n$, let $a_{i}$ denote the unique minimal lattice
  generator of $\rho_{i}$ in $\Sigma$ and let $\ell_{i}$ be the
  positive integer satisfying the relation $\bar b_{i} = \ell_{i}
  a_{i}$.  The Jurkiewicz-Danilov Theorem (see Page~134 in \cite{Oda})
  states that there is a surjective homomorphism of graded rings from
  $\mathbb{Q}[x_{1}, \dotsc, x_{n}]$ to $A^{*}\big( X(\Sigma) \big)$
  given by $x_{i} \mapsto D_{i}$ where $D_{i}$ is the torus invariant
  Weil divisor on $X(\Sigma)$ associated with $\rho_{i}$.  The kernel
  of this map is the ideal $I_{\Sigma}$ plus the ideal generated by
  the linear relations $\sum_{i=1}^{n} \theta(a_{i}) \cdot x_{i}$ for
  all $\theta \in {N}^{\star}$.  Example~6.7 in \cite{Vistoli}
  establishes a natural isomorphism
  $A^{*}\big(\mathcal{X}(\bm{\Sigma})\big) \cong A^{*}\big( X(\Sigma)
  \big)$ defined by $c_{1}(L_{i}) \mapsto \ell_{i}^{-1} \cdot D_{i}$.
  Since we have $\sum_{i=1}^{n} \theta(a_{i}) \cdot \ell_{i} \cdot
  x_{i} = \sum_{i=1}^{n} \theta(\bar b_{i}) \cdot x_{i}$ for all
  $\theta \in {N}^{\star}$, the composition of these two isomorphism
  establishes the claim.
\end{proof}

This lemma allow us to establish Theorem~\ref{thm:main} at the level
of $\mathbb{Q}$-graded $\mathbb{Q}$-vector spaces.  More precisely, we
prove the following the result.  If $M$ is a $\mathbb{Q}$-graded
module and $c$ is a rational number, then we write $M[c]$ for the
$c$-th shift of $M$; its defined by the formula $M[c]_{c'} = M_{c' +
  c}$.

\begin{proposition} \label{pro:moduleIso}
  If $\mathcal{X}(\bm{\Sigma})$ is a complete toric Deligne-Mumford
  stack, then there is an isomorphism of $\mathbb{Q}$-graded
  $\mathbb{Q}$-vector spaces:
  \[
  \frac{\mathbb{Q}[N]^{\bm{\Sigma}}}{\left\langle \sum_{i=1}^{n}
      \theta(b_{i}) \cdot y^{b_{i}} : \theta \in N^{\star}
    \right\rangle} \cong \bigoplus_{v \in \BBox(\bm{\Sigma})}
  A^{*}\big( \mathcal{X}( \bm{\Sigma / \sigma}(\bar v)) \big) \bigl[
  \deg(y^{v}) \bigr] \, .
  \]
\end{proposition}

\begin{proof}
  The definition of $S_{\bm{\Sigma}}$ and $\BBox(\bm{\Sigma})$ implies
  that $\mathbb{Q}[N]^{\bm{\Sigma}} = \bigoplus_{v \in
    \BBox(\bm{\Sigma})} y^{v} \cdot S_{\bm{\Sigma}}$.  We first
  analyze the individual summands.  Fix an element $v \in
  \BBox(\bm{\Sigma})$ and let $\tau := \sigma(\bar{v})$ be the minimal
  cone in $\Sigma$ containing $\bar{v}$.  It follows from the
  definition of multiplication in the deformed group ring that $y^{v}
  \cdot S_{\bm{\Sigma}}$ is isomorphic to the quotient of
  $S_{\bm{\Sigma}}$ by the ideal generated by the elements $y^{c}$
  where $c$ lies outside the cones in $\Sigma$ containing $\tau$.
  
  Let $S_{\bm{\Sigma/\tau}}$ denote the subring of
  $\mathbb{Q}[N(\tau)]^{\bm{\Sigma/\tau}}$ generated by
  $y^{\tilde{b}_i}$ for $\rho_i\in\link(\tau)$.  By renumbering the
  rays in $\Sigma$, we may assume that $\tilde{\rho}_{1}, \dotsc,
  \tilde{\rho}_{\ell}$ are the rays in $\link(\tau)$.  Recall that
  $\tilde{b}_{i}$ is the image of $b_{i}$ in $N(\tau)$.  For each ray
  $\rho_{i} \in \tau$, choose an element $\theta_{i} \in N^{\star}$
  such that $\theta_{i}(b_{i}) = 1$ and $\theta_{i}(b_{j}) = 0$ for
  all $\bar{b}_{i} \neq \bar{b}_{j} \in \tau$.  Consider the map
  defined by
  \[
  y^{b_{i}} \mapsto
  \begin{cases} 
    y^{\tilde{b}_{i}} & \text{for $\rho_{i} \subseteq \link(\tau)$;}\\
    - \sum_{j=1}^{\ell} \theta_{i}(b_{j}) \cdot y^{\tilde{b}_{j}} &
    \text{for $\rho_{i} \subseteq \tau$;}\\
    0 & \text{for $\rho_{i} \not\subseteq \tau \cup \link(\tau)$.}
  \end{cases}
  \]
  Since this map is compatible with the multiplicative structures on
  $S_{\bm{\Sigma}}$ and $S_{\bm{\Sigma / \tau}}$, it induces a
  surjective homomorphism from $S_{\bm{\Sigma}}$ to $S_{\bm{\Sigma /
      \tau}}$.  Clearly, the kernel contains the elements
  $\theta_{i}(b_{i}) \cdot y^{b_{i}} + \sum_{j = 1}^{\ell}
  \theta_{i}(b_{j}) y^{b_{j}}$ for all $\rho_{i} \in \tau$ and the
  elements $y^{c}$ where $c$ lies outside the cones in $\Sigma$
  containing $\tau$.  Given any other element of the kernel, we can
  use these relations to obtain a linear combination of monomials
  $y^{w}$ with $\bar{w} \in \link(\tau)$ which also belongs to the
  kernel.  However, this is only possible if all the coefficients of
  $y^{w}$ are zero which implies that the given elements generate the
  kernel.
  
  Since Lemma~\ref{lem:basicChow} establishes that
  \[
  \frac{S_{\bm{\Sigma/\tau}}} {\langle \sum_{i = 1}^{\ell}
    \tilde{\theta} ( \tilde{b}_{i}) \cdot y^{\tilde{b}_{i}} :
    \tilde{\theta} \in N(\tau)^{\star} \rangle} \cong A^{*}\big(
  \mathcal{X}( \bm{\Sigma / \tau}) \big) \, ,
  \]
  we have a surjective $\mathbb{Q}$-graded $\mathbb{Q}$-linear map
  from $y^{v} \cdot S_{\bm{\Sigma}}$ to $A^{*}\big( \mathcal{X}(
  \bm{\Sigma / \tau}) \big)[\deg(y^{v})]$ whose kernel is generated by
  the elements $\theta_{i}(b_{i}) \cdot y^{b_{i}} + \sum_{j =
    1}^{\ell} \theta_{i}(b_{j}) y^{b_{j}}$ for all $\rho_{i} \in \tau$
  and the pullbacks of the linear relations $\sum_{i = 1}^{\ell}
  \tilde{\theta} ( \tilde{b}_{i}) \cdot y^{\tilde{b}_{i}}$ where
  $\tilde{\theta} \in N(\tau)^{\star}$.  Finally, taking the direct
  sum over all $v \in \BBox(\bm{\Sigma})$ produces a surjective
  $\mathbb{Q}$-graded $\mathbb{Q}$-linear map from
  $\mathbb{Q}[N]^{\bm{\Sigma}}$ to $\bigoplus_{v \in
    \BBox(\bm{\Sigma})} A^{*}\big( \mathcal{X}( \bm{\Sigma /
    \sigma}(\bar v)) \big) \bigl[ \deg(y^{v}) \bigr]$ whose kernel is
  generated by the elements $\sum_{i = 1}^{n} \theta(b_{i}) \cdot
  y^{b_{i}}$ where $\theta \in N^{\star}$.
 \end{proof}

\begin{remark}
  Although the elements $\theta_{i}$ in the proof of
  Proposition~\ref{pro:moduleIso} are not uniquely determined, the
  possible choices differ by elements in $N(\tau)^{\star}$.  It
  follows that the surjection from $y^{v} \cdot S_{\bm{\Sigma}}$ to
  $S_{\bm{\Sigma / \tau}}[\deg(y^{v})]$ is not canonically defined,
  but the surjection from $y^{v} \cdot S_{\bm{\Sigma}}$ to $A^{*}\big(
  \mathcal{X}( \bm{\Sigma / \tau}) \big) \bigl[ \deg(y^{v}) \bigr]$
  is.
\end{remark}

\begin{remark} \label{rem:grading}
  The degree shift in Proposition~\ref{pro:moduleIso} is also called
  the \emph{age} of the component of the inertia stack; Subsection~7.1
  in \cite{AGV}.
\end{remark}

\section{The Product Structure on $A_{orb}^{*} \bigl(
  \mathcal{X}(\bm{\Sigma}) \bigr)$} \label{sec:orbchow}
\setcounter{equation}{0}

In this section, we study multiplication in $A_{orb}^{*} \bigl(
\mathcal{X}(\bm{\Sigma}) \bigr)$.  Specifically, we complete the proof
of Theorem~\ref{thm:main} by showing that multiplication in the
deformed group ring coincides with the orbifold product.

To compare the two products, we first give a combinatorial description
of the moduli space $\mathcal{K} := \mathcal{K}_{0,3} \bigl(
\mathcal{X} (\bm{\Sigma}), 0 \bigr)$ of $3$-pointed twisted stable
maps of genus zero and degree zero to $\mathcal{X}(\bm{\Sigma})$.  The
moduli space $\mathcal{K}$ is a smooth proper Deligne-Mumford stack
with a projective coarse moduli space; see Theorem~3.6.2 in
\cite{AGV}.  In addition, Lemma~6.2.1 in \cite{AGV} gives three
evaluation maps denoted $\ev_{i} \colon {\mathcal{K}} \rightarrow
{\mathcal{I}} \bigl( \mathcal{X}(\bm{\Sigma}) \bigr)$ for $1 \leq i
\leq 3$.  Proposition~\ref{pro:inertia} shows that $\mathcal{I} \bigl(
\mathcal{X} (\bm{\Sigma}) \bigr) = \coprod_{v \in \BBox(\bm{\Sigma})}
\mathcal{X} \bigl( \bm{\Sigma / \sigma}(\bar{v}) \bigr)$, so we can
index the components of $\mathcal{K}$ by the images of the evaluation
maps.  Let $\mathcal{K}_{v_{1}, v_{2}, v_{3}}$ be the component of
$\mathcal{K}$ such that $\ev_{i}$ maps to $\mathcal{X} \bigl(
\bm{\Sigma / \sigma}(\bar{v}_{i}) \bigr)$ for $1 \leq i \leq 3$.

For brevity, we write $v_{1} + v_{2} + v_{3} \equiv 0$ to indicate
that there exists a cone $\sigma \in \Sigma$ containing $\bar{v}_{i}$
for $1 \leq i \leq 3$ such that the sum $v_{1} + v_{2} + v_{3}$
belongs to the subgroup $N_{\sigma}$ in $N$.

\begin{proposition} \label{pro:K_0,3}
  If $\mathcal{X}(\bm{\Sigma})$ is a complete toric Deligne-Mumford
  stack, then
  \[
  \mathcal{K} = \coprod_{\stackrel{(v_{1}, v_{2}, v_{3}) \in
      \BBox(\bm{\Sigma})^3}{v_{1} + v_{2} + v_{3} \equiv 0}}
  \mathcal{X} \bigl( \bm{\Sigma / \sigma}(\bar{v}_{1}, \bar{v}_{2},
  \bar{v}_{3}) \bigr) \, ,
  \]
  where $\sigma(\bar{v}_{1}, \bar{v}_{2}, \bar{v}_{3})$ is the minimal
  cone in $\Sigma$ containing $\bar{v}_{1}$, $\bar{v}_{2}$ and
  $\bar{v}_{3}$.
\end{proposition}

\begin{proof}
  We begin by examining the geometric points of $\mathcal{K}$.  A
  $\mathbb{C}$-valued point of $\mathcal{K}$ is a representable
  morphism $f$ from a twisted curve $\mathcal{C}$ to
  $\mathcal{X}(\bm{\Sigma})$ such that the induced map on coarse
  moduli spaces sends $\mathbb{P}^{1}$ to a point $x \in X(\Sigma)$.
  Hence, the map $f$ factors through a closed substack $\mathcal{B}G'$
  in $\mathcal{X}(\bm{\Sigma})$ where $G' \subseteq G$ is the isotropy
  group of $x \in \mathcal{X}(\bm{\Sigma})$ and $\mathcal{B} G'$ is
  the classifying stack $[x/G']$.  Corollary~1.6.2 in
  \cite{LaumonMoret} shows that the morphism from $\mathcal{C}$ to
  $\mathcal{B}G'$ is also representable which implies that the fibered
  product $\widehat{C} := \mathcal{C} \times_{\mathcal{B}G'} x$ is a
  scheme.  Since $\mathcal{C}$ is smooth, we see that $\widehat{C}$ is
  a smooth curve, although it is typically disconnected.  Let $H$ be
  the subgroup of $G'$ that acts trivially on the set of connected
  components of $\widehat{C}$.  Since $G'$ is abelian, the group $H$
  is the stabilizer of each connected component of $\widehat{C}$.  By
  choosing a connected component $C$ of $\widehat{C}$, we obtain
  $\mathcal{C} \cong [C/H]$.  Assuming the points $\{ 0, 1, \infty \}$
  in $\mathbb{P}^{1}$ correspond to the markings on $\mathcal{C}$, the
  properties of a twisted curve imply that the map $\mathcal{C}
  \rightarrow \mathbb{P}^{1}$ is an isomorphism over $\mathbb{P}^{1} -
  \{0,1, \infty\}$.  It follows that $C$ is a proper smooth Galois
  cover of $\mathbb{P}^{1}$ with Galois group $H$ branched over $0$,
  $1$ and $\infty$.  Specifically, if $\gamma_{1}, \gamma_{2},
  \gamma_{3}$ are the generators of the fundamental group of
  $\mathbb{P}^{1} - \{ 0, 1, \infty \}$ corresponding to
  counterclockwise loops around $0$, $1$, $\infty$ respectively, then
  $C$ is induced by a homomorphism $\pi_{1}(\mathbb{P}^{1} - \{ 0, 1,
  \infty \}) \rightarrow G$ sending $\gamma_{1}$ to $g_{i}$ such that
  $g_{1} \cdot g_{2} \cdot g_{3} = 1$ and $g_{i}$ generate $H$ as a
  subgroup of $G$.
  
  By definition, the map $\ev_{i}$ is induced by the representable
  morphism from the cyclotomic gerbe in $\mathcal{C}$ lying over the
  corresponding point in $\mathbb{P}^{1}$ to
  $\mathcal{X}(\bm{\Sigma})$; recall that over $\mathbb{C}$ the
  inertia stack $\mathcal{I} \bigl( \mathcal{X}(\bm{\Sigma}) \bigr)$
  is canonically isomorphic to the stack of representable morphisms
  from a constant cyclotomic gerbe to $\mathcal{X}(\bm{\Sigma})$.
  Hence, the evaluation map $\ev_{i}$ sends $f$ to the geometric point
  $(x, g_{i})$ in the inertia stack.  Because $g_{i}$ belongs to the
  isotropy group of $x$, it fixes a point in $Z$.  Thus,
  Lemma~\ref{lem:boxFixed} shows that $g_{i}$ corresponds to an
  element $v_{i} \in \BBox(\bm{\Sigma})$ and $\ev_{i}$ maps to the
  component $[Z^{g_{i}} / G] = \mathcal{X} \bigl( \bm{\Sigma /
    \sigma}(\bar{v}_{i}) \bigr)$ of the inertia stack.  Moreover, the
  condition that $g_{1} \cdot g_{2} \cdot g_{3} = 1$ means that there
  exists a cone $\sigma \in \Sigma$ containing $\bar{v}_{1},
  \bar{v}_{2}, \bar{v}_{3}$ and the sum $v_{1} + v_{2} + v_{3}$
  belongs to the subgroup $N_{\sigma}$ in $N$.  Therefore, the
  component $\mathcal{K}_{v_{1}, v_{2}, v_{3}}$ is nonempty if and
  only if $v_{1} + v_{2} + v_{3} \equiv 0$.
  
  The morphisms $\ev_{i} \colon \mathcal{K}_{v_{1}, v_{2}, v_{3}}
  \rightarrow \mathcal{X} \bigl( \bm{\Sigma / \sigma}(\bar{v}_{i})
  \bigr)$ are compatible with the inclusion maps into
  $\mathcal{X}(\bm{\Sigma})$ for $1 \leq i \leq 2$ which yields a
  morphism
  \[
  e \colon \mathcal{K}_{v_{1}, v_{2}, v_{3}} \rightarrow \mathcal{X}
  \bigl( \bm{\Sigma / \sigma}(\bar v_{1}) \bigr)
  \times_{\mathcal{X}(\bm{\Sigma})} \mathcal{X} \bigl( \bm{\Sigma /
    \sigma}(\bar{v}_{2}) \bigr) = [Z^{g_{1}} / G] \times_{[Z/G]}
  [Z^{g_{2}} / G] \quad .
  \]
  Because $H$ is the subgroup of $G$ generated by $g_{1}$ and $g_{2}$
  (note: $g_{3} = g_{1}^{-1} g_{2}^{-1}$), we have $Z^{g_{1}}
  \times_{Z} Z^{g_{2}} = Z^{\langle g_{1}, g_{2} \rangle} = Z^{H}$.
  It follows that $[Z^{g_{1}} / G] \times_{[Z/G]} [Z^{g_{2}} / G] =
  [Z^{H}/ G]$.  Our analysis of the geometric points of $\mathcal{K}$
  shows that $e$ induces a bijection between the
  $\mathbb{C}$\nobreakdash-valued points of the coarse moduli spaces
  of $\mathcal{K}_{v_{1}, v_{2}, v_{3}}$ and $[Z^{H} / G]$.  Since
  both $\mathcal{K}_{v_{1},v_{2},v_{3}}$ and $[Z^{H} / G]$ are smooth
  Deligne-Mumford stacks, their coarse moduli spaces have at worst
  quotient singularities.  Applying Theorem~VI.1.5 in \cite{Kollar},
  we deduce that, in fact, $e$ produces an isomorphism between the
  coarse moduli spaces.
  
  To prove that $e$ is an isomorphism of stacks, it remains to show
  that $e$ gives an isomorphism between the isotropy groups of
  $\mathbb{C}$-valued points.  Indeed, since $\mathcal{K}$ is smooth
  (see page~18 in \cite{AGV}) and $e$ is representable, the
  isomorphism follows from a similar statement for the lifting of $e$
  to the atlases.  Proposition~7.1.1 in \cite{ACV} indicates that the
  automorphism group of a twisted stable curve is the direct product
  of the automorphism groups of the nodes which implies that our curve
  $\mathcal{C}$ has only the trivial automorphism.  Hence, an isotropy
  of the twisted stable map $f \colon \mathcal{C} \rightarrow
  \mathcal{B}G' \subseteq \mathcal{X}(\bm{\Sigma})$ corresponds to a
  diagram
  \[
  \begin{CD}
    E @>{\phi}>> E'\\
    @VVV @VVV\\
    \mathcal{C} @= \mathcal{C} 
  \end{CD}
  \]
  where $\phi$ is a $G'$-equivariant map of principal $G'$-bundles
  over $\mathcal{C}$.  Since $\mathcal{C}$ is connected, the map
  $\phi$ is multiplication by an element of $G'$.  Therefore, the
  isotropy group of the map $f$ is precisely $G'$ which completes the
  proof.
\end{proof}

Proposition~\ref{pro:K_0,3} also provides a presentation for the
universal twisted stable curve over $\mathcal{K}$.  To describe the
universal curve, we focus on the component $\mathcal{K}_{v_{1}, v_{2},
  v_{3}}$.  As above, we write $H$ for the subgroup of $G$
corresponding to $\{ v_{1}, v_{2}, v_{3} \}$ and $C \rightarrow
\mathbb{P}^{1}$ for the associated Galois cover.  Consider the
quotient stack
\[
\mathcal{U}_{v_{1}, v_{2}, v_{3}} := [(Z^{H} \times C) / (G \times H)]
= [Z^{H} / G] \times [C / H] \quad .
\]  
If $S$ is a scheme, then the objects in $\mathcal{U}_{v_{1}, v_{2},
  v_{3}}(S)$ are principal $(G \times H)$-bundles $E \rightarrow S$
with a $(G \times H)$-equivariant map $E \rightarrow Z^{H} \times C$.
The \emph{twisted projection map} $\pi$ from $\mathcal{U}_{v_{1},
  v_{2}, v_{3}}$ to $\mathcal{K}_{v_{1}, v_{2}, v_{3}} = [Z^{H} / G]$
is defined as follows: If $H$ acts on $E$ via the map $h \mapsto
(h^{-1}, h) \in G \times H$, then $E/H$ is a principal $G$-bundle over
$S$.  To obtain an object in $\mathcal{K}_{v_{1}, v_{2}, v_{3}}(S)$,
observe that the $(G \times H)$-equivariant map $E \rightarrow Z^{H}
\times C$ induces a $G$\nobreakdash-equivariant map from $E/H$ to
$Z^{H}$.  By verifying that $\pi$ is compatible with morphisms in
$\mathcal{U}_{v_{1}, v_{2}, v_{3}}(S)$ and $\mathcal{K}_{v_{1}, v_{2},
  v_{3}}(S)$, we conclude that $\pi$ is a morphism of stacks.  With
these definitions, we have

\begin{corollary} \label{cor:universal}  
  The universal twisted stable curve over $\mathcal{K}_{v_{1}, v_{2},
    v_{3}} \cong [Z^{H} / G]$ is given by the twisted projection map
  $\pi \colon \mathcal{U}_{v_{1}, v_{2}, v_{3}} = [(Z^{H} \times C)/(G
  \times H)] \rightarrow [Z^{H}/G]$.
\end{corollary}

\begin{proof}
  Fix a map $S \rightarrow [Z^{H} / G]$ where $S$ is a scheme and
  consider the fibered product $\mathcal{D} := \mathcal{U}_{v_{1},
    v_{2}, v_{3}} \times_{[Z^{H} /G]} S$.  Assuming that $S
  \rightarrow [Z^{H} / G]$ corresponds to the principal
  $G$\nobreakdash-bundle $E \rightarrow S$ with a $G$-equivariant map
  $E \rightarrow Z^{H}$, we have $\mathcal{D} = [ (E \times C) / (G
  \times H)]$ where the $(G \times H)$-action is given by $(e,c,g,h)
  \mapsto (e \cdot gh^{-1}, c \cdot h)$.  The twisted projection map
  $\pi$ induces a map $\mathcal{D} \rightarrow [E/G] = S$.  Because
  the \emph{anti-diagonal} action of $H$ on $E \times C$ is free, the
  quotient $Y := (E \times C)/H$ is a scheme.  Hence, we have
  $\mathcal{D} = [Y/G]$ where the $G$-action on $Y$ is induced by the
  action on $E \times C$.  Since $H$ acts trivially on $Z^{H}$, the
  $G$\nobreakdash-equivariant map $E \rightarrow Z^{H}$ induces a
  $G$-equivariant map $Y \rightarrow Z^{H}$ which shows that
  $\mathcal{D}$ maps to $[Z^{H}/ G] \subseteq [Z / G] =
  \mathcal{X}(\bm{\Sigma})$.  Moreover, if $R = R_{1} + R_{2} + R_{3}$
  is the ramification divisor of the Galois cover $C \rightarrow
  \mathbb{P}^{1}$, then the image of the open set $E \times (C-R)$
  gives an open substack of $[Y/G]$ which is isomorphic to $S \times
  \bigl( \mathbb{P}^{1} - \{ 0, 1, \infty \} \bigr)$.  By definition,
  the evaluation map $\ev_{i}$ from $\mathcal{D}$ to the inertia stack
  $\mathcal{I} \bigl( \mathcal{X}(\bm{\Sigma}) \bigr)$ arises from the
  representable morphism from $[(E \times R_{i}) / (G \times H)]$ to
  $\mathcal{X}(\bm{\Sigma})$.  In particular, $\ev_{i}$ is induced by
  the closed embedding $[Z^{H}/ G] \rightarrow [Z^{g_{i}} / G] \cong
  \mathcal{X} \bigl( \bm{\Sigma / \sigma}(\bar{v}_{i}) \bigr)$.  We
  conclude that $\mathcal{U}_{v_{1}, v_{2}, v_{3}}$ is a family of
  twisted stable curves over $[Z^{H} / G]$ with a map $f \colon
  \mathcal{U}_{v_{1}, v_{2}, v_{3}} \rightarrow
  \mathcal{X}(\bm{\Sigma})$ and evaluation maps $\ev_{i} \colon
  \mathcal{U}_{v_{1}, v_{2}, v_{3}} \rightarrow \mathcal{X} \bigl(
  \bm{\Sigma / \sigma}(\bar{v}_{i}) \bigr) \subseteq \mathcal{I}
  \bigl( \mathcal{X}(\bm{\Sigma}) \bigr)$ for $1 \leq i \leq 3$.
  
  Let $\mathcal{U}'$ denote the universal family of twisted stable
  curves over $\mathcal{K}_{v_{1}, v_{2}, v_{3}}$.  By the universal
  mapping property of $\mathcal{U}'$, there exists a map $\mu \colon
  [Z^{H} / G] \rightarrow \mathcal{K}_{v_{1}, v_{2}, v_{3}}$ such that
  \[
  \begin{CD}
    \mu^{*}(\mathcal{U}') @>>> \mathcal{U}' \\
    @VVV @VVV \\
    [Z^{H} / G] @>{\mu}>> \mathcal{K}_{v_{1}, v_{2}, v_{3}}
  \end{CD}
  \]
  is a Cartesian diagram.  Combining definition of $e$ with the first
  paragraph, we see that $e \circ \mu = \text{id}$. Since
  Proposition~\ref{pro:K_0,3} shows that $e$ is an isomorphism, we
  conclude that $\mu$ is also an isomorphism and $\mathcal{U}_{v_{1},
    v_{2}, v_{3}}$ is isomorphic to $\mathcal{U}'$.
\end{proof}

Next, we describe the virtual fundamental class on $\mathcal{K}$.
Recall that $L_{k}$ denotes the line bundle on
$\mathcal{X}(\bm{\Sigma})$ corresponding to the line bundle
$\mathbb{C} \times Z$ on $Z$ where the $G$-action on $\mathbb{C}$
given by the $k$-th component $\alpha_{k}$ of $\alpha \colon G
\rightarrow (\mathbb{C}^{*})^{n}$.

\begin{proposition}\label{pro:virtual}
  Let $\mathcal{K}_{v_{1}, v_{2}, v_{3}}$ be a component of the moduli
  space $\mathcal{K}$.  If the integers $m_{k} \in \{ 1, 2 \}$ are
  defined by the relation $v_{1} + v_{2} + v_{3} = \sum_{\rho_{i} \in
    \sigma(\bar{v}_{1}, \bar{v}_{2}, \bar{v}_{3})} m_{k} b_{k}$ in
  $N$, then the virtual fundamental class of the component
  $\mathcal{K}_{v_{1}, v_{2}, v_{3}}$ is
  \[
  \prod_{m_k = 2} c_{1}(L_{k}) \big|_{\mathcal{X} \left( \bm{\Sigma /
        \sigma}(\bar{v}_{1}, \bar{v}_{2}, \bar{v}_{3}) \right)} \quad
  .
  \]
\end{proposition}

\begin{proof}
  Let $f$ be the natural map from $\mathcal{U}_{v_{1}, v_{2}, v_{3}}$
  to $\mathcal{X}(\bm{\Sigma})$ and let $\pi \colon
  \mathcal{U}_{v_{1}, v_{2}, v_{3}} \rightarrow [Z^{H} / G]$ be the
  twisted projection map.  Since $\mathcal{K}_{v_{1}, v_{2}, v_{3}}$
  is smooth, the virtual fundamental class of $\mathcal{K}$ is given
  by the top Chern class of the bundle $R^{1} \pi_{*} f^{*}
  (T_{\mathcal{X}(\bm{\Sigma})})$; see Section~6.2 in \cite{AGV}.  To
  calculate this Chern class, observe that the pullback of the tangent
  bundle $f^{*} (T_{\mathcal{X}(\bm{\Sigma})})$ corresponds to a $(G
  \times H)$-equivariant bundle $\mathcal{V}$ on $Z^{H} \times C$;
  $\mathcal{V}$ is a trivial vector bundle of rank $n$ where the $(G
  \times H)$-action is induced by the map $\alpha \colon G \rightarrow
  (\mathbb{C}^{*})^{n}$ on it basis.  Let $p \colon Z^{H} \times C
  \rightarrow Z^{H}$ be the projection map and let $p_{*}^{H}$ be the
  invariant pushforward (pushing forward and taking invariant
  sections).  Since the associated derived functor $R^{1}p_{*}^{H}$
  sends $(G \times H)$-equivariant sheaves on $Z^{H} \times C$ to
  $G$-equivariant sheaves on $Z^{H}$, it suffices to compute $R^{1}
  p_{*}^{H} (\mathcal{V})$.
  
  Let $\mathcal{W}_{k}$ be the trivial line bundle on $Z^{H} \times C$
  with $(G \times H)$-action induced by the $k$-th component
  $\alpha_{k}$ of $\alpha \colon G \rightarrow (\mathbb{C}^{*})^{n}$
  and consider the following exact sequence of vector bundles on
  $Z^{H} \times C$:
  \[
  0 \longrightarrow p^{*}( T_{Z^{H}}) \longrightarrow \mathcal{V}
  \longrightarrow \bigoplus_{\rho_{k} \in \sigma (\bar{v}_{1},
    \bar{v}_{2}, \bar{v}_{3})} \mathcal{W}_{k} \longrightarrow 0 \quad
  .
  \]
  Since the $H$-invariant part of $R^{1}p_{*} p^{*}(T_{Z^{H}}) = R^{1}
  p_{*}(\mathcal{O}_{Z^{H} \times C}) \otimes T_{Z^{H}}$ is trivial,
  it suffices to calculate $R^{1}p_{*}^{H}(\mathcal{W}_{k})$.  Given a
  point $z \in Z^{H}$, the restriction of $\mathcal{W}_{k}$ to $z
  \times C$ is isomorphic to the trivial line bundle $\mathcal{L}_{k}$
  on $C$ with the the $H$-action induced by $\alpha_{k}$.  Since the
  Leray spectral sequence degenerates, we have $H^{1}(C,
  \mathcal{L}_{k}) \cong H^{1}\bigl( \mathbb{P}^{1}, (p')_{*}^{H}
  (\mathcal{L}_{k}) \bigr)$ where $p' \colon C \rightarrow
  \mathbb{P}^{1}$ is the Galois cover.  Because $v_{j} \in
  \BBox(\bm{\Sigma})$ for $1 \leq j \leq 3$, there are $a_{j,k} \in
  \mathbb{Q}$ such that $0 \leq a_{j,k} \leq 1$ and $\bar{v}_{j} =
  \sum a_{j,k} \bar{b}_{k}$ where $\rho_{k} \in \sigma(\bar{v}_{1},
  \bar{v}_{2}, \bar{v}_{3})$.  By hypothesis, we have $v_{1} + v_{2} +
  v_{3} \equiv 0$ which means that $a_{1,k} + a_{2,k} + a_{3,k}$ is an
  integer between $0$ and $2$.  Lemma~\ref{lem:boxFixed} establishes
  that $v_{j}$ corresponds to an element $g_{j} \in G$ and the proof
  of Proposition~\ref{pro:open} shows that $\alpha_{k}(g_{j}) = \exp(2
  \pi \sqrt{-1} a_{j,k})$.  From this explicit description of the
  $H$-action on $\mathcal{L}_{k}$, it follows that
  $(p')_{*}^{H}(\mathcal{L}_{k}) \cong p_{*}^{H} \bigl(
  \mathcal{W}_{k} \big|_{z \times C} \bigr)$ is isomorphic to
  $\mathcal{O}_{\mathbb{P}^{1}}(- a_{1,k} - a_{2,k} - a_{3,k})$.
  Since
  \[
  \dim H^{1} \bigl( \mathbb{P}^{1}, \mathcal{O}_{\mathbb{P}^{1}}(-
  a_{1,k} - a_{2,k} - a_{3,k}) \bigr) = 1 \qquad \text{when $a_{1,k} +
    a_{2,k} + a_{3,k} = 2$,}
  \] 
  we deduce that $R^{1}p_{*}^{H}(\mathcal{W}_{k})$ is the line bundle
  $\mathbb{C} \times Z$ on $Z$ where the $G$-action on $\mathbb{C}$
  given by the $k$-th component $\alpha_{k}$.  When $a_{1,k} + a_{2,k}
  + a_{3,k} \neq 2$, the cohomology group vanishes and
  $R^{1}p_{*}^{H}(\mathcal{W}_{k})$ is zero.  Therefore, we have
  \[
  R^{1} \pi_{*} f^{*} (T_{\mathcal{X}(\bm{\Sigma})}) \cong
  \bigoplus_{m_k = 2} L_{k} \big|_{[Z^{H}/G]} 
  \]
  and taking the top Chern class completes the proof.
\end{proof}

\begin{remark}  
  The virtual classes in Proposition~\ref{pro:virtual} are analogous
  to the classes $c(g,h)$ in \cite{FantechiGottsche}.  However, we do
  not need to use the language of parabolic bundles because we give an
  explicit description for the $H$-action on the trivial line bundles
  $\mathcal{L}_{k}$.
\end{remark}

We end this section with a proof of Theorem~\ref{thm:main}.  Let
$\iota \colon \mathcal{I} \bigl( \mathcal{X}(\bm{\Sigma}) \bigr)
\rightarrow \mathcal{I} \bigl( \mathcal{X}(\bm{\Sigma}) \bigr)$ denote
the natural involution on the inertia stack defined by $(x, \phi)
\mapsto (x, \phi^{-1})$ and let $\check{\ev}_{3} := \iota \circ \,
\ev_{3}$ be the twisted evaluation map; see Section~4.5 in \cite{AGV}.
If $\gamma_{1}, \gamma_{2} \in A^{*}\bigl( \mathcal{I} \bigl(
\mathcal{X}(\bm{\Sigma}) \bigr)$, then the orbifold product
(Definition~6.2.2 in \cite{AGV}) is 
\[
\gamma_{1} * \gamma_{2} := (\check{\ev}_{3})_{*} \bigl( \ev_{1}^{*}
(\gamma_{1}) \cap \ev_{2}^{*} (\gamma_{2}) \cap [ \mathcal{K} ]^{vir}
\bigr)
\] 
where $[\mathcal{K}]^{vir}$ denotes the virtual fundamental class on
$\mathcal{K}$.  This definition agrees with the definition of the
quantum product in degree zero.

\begin{remark} \label{rem:nonempty}
  Proposition~\ref{pro:K_0,3} shows that the component
  $\mathcal{K}_{v_{1}, v_{2}, v_{3}}$ of the moduli stack is nonempty
  if and only if $v_{1} + v_{2} + v_{3} \equiv 0$.  Hence, if
  $\gamma_{1} \in A^{*} \bigl( \mathcal{X}( \bm{ \Sigma /
    \sigma}(\bar{v}_{1}) \bigr)$ and $\gamma_{2} \in A^{*} \bigl(
  \mathcal{X}( \bm{ \Sigma / \sigma}(\bar{v}_{2}) \bigr)$, then the
  orbifold product $\gamma_{1} * \gamma_{2}$ is nonzero only if there
  is a cone in $\Sigma$ containing $\bar{v}_{1}$ and $\bar{v}_{2}$.
\end{remark}

\begin{proof}[Proof of Theorem~\ref{thm:main}]
  By combining Proposition~\ref{pro:inertia} and
  Proposition~\ref{pro:moduleIso}, we obtain the following isomorphism
  of $\mathbb{Q}$-graded $\mathbb{Q}$-vector spaces:
  \[
  A_{orb}^{*} \bigl( \mathcal{X}(\bm{\Sigma}) \bigr) = \bigoplus_{v
    \in \BBox(\bm{\Sigma})} A^{*}\big( \mathcal{X}( \bm{\Sigma /
    \sigma}(\bar v)) \big) \bigl[ \deg(y^{v}) \bigr] \cong
  \frac{\mathbb{Q}[N]^{\bm{\Sigma}}}{\left\langle \sum_{i=1}^{n}
      \theta(b_{i}) \cdot y^{b_{i}} : \theta \in N^{\star}
    \right\rangle} \, .
  \]
  It remains to show that the orbifold product agrees with the product
  structure on the deformed group ring.  Since the elements of
  $\BBox(\bm{\Sigma})$ generates $ A_{orb}^{*} \bigl(
  \mathcal{X}(\bm{\Sigma}) \bigr)$ as a module over the $y^{b_{i}}$,
  it suffices to show that $y^{c} * y^{b_{i}} = y^{c} \cdot y^{b_{i}}$
  and $y^{v_{1}} * y^{v_{2}} = y^{v_{1}} \cdot y^{v_{2}}$ where $c \in
  N$ and $v_{1}, v_{2} \in \BBox(\bm{\Sigma})$.
  
  We first consider the product $y^{c} * y^{b_{i}}$ where $c \in N$.
  By taking advantage of the linear relations $\sum_{i=1}^{n}
  \theta(b_{i}) \cdot y^{b_{i}}$ for $\theta \in N^{\star}$, we reduce
  to the case that $b_{i}$ does not lie in the minimal cone
  $\sigma(\bar{c})$ containing $\bar{c}$.  Let $v$ be the
  representative of $c$ in $\BBox(\bm{\Sigma})$.  By
  Remark~\ref{rem:nonempty}, the only contribution to the product
  $y^{c} * y^{b_{i}}$ comes from the component $\mathcal{K}_{v,0,v'}$
  where $v' \in \BBox(\bm{\Sigma})$ is defined by the equation $v + v'
  = \sum_{\rho_{i} \in \sigma(\bar{c})} b_{i}$.  Hence,
  $\mathcal{K}_{v,0,v'}$ is isomorphic to $\mathcal{X} \bigl(
  \bm{\Sigma / \sigma}(\bar{c}) \bigr)$, both $\ev_{1},
  \check{\ev}_{3} \colon \mathcal{X} \bigl( \bm{\Sigma /
    \sigma}(\bar{c}) \bigr) \rightarrow \mathcal{X} \bigl( \bm{\Sigma
    / \sigma}(\bar{c}) \bigr)$ are the identity map and $\ev_{2}
  \colon \mathcal{X} \bigl( \bm{\Sigma / \sigma}(\bar{c}) \bigr)
  \rightarrow \mathcal{X} ( \bm{\Sigma})$ is the closed embedding.
  The restriction of $y^{b_{i}}$ from $\mathcal{X}(\bm{\Sigma})$ to
  $\mathcal{X} \bigl( \bm{\Sigma / \sigma}(\bar{c}) \bigr)$ is equal
  to $y^{\tilde{b}_{i}}$ if $\bar{b}_{i}$ and $\sigma(\bar{c})$ lie in
  a cone of $\Sigma$ and is equal to zero otherwise.  Since
  Proposition~\ref{pro:virtual} shows the the virtual fundamental
  class is $1$, if $\ev_{2}^{*}(y^{b_{i}}) \neq 0$ then $y^{c} *
  y^{b_{i}}$ is simply multiplication in $A^{*}(\mathcal{X} \bigl(
  \bm{\Sigma / \sigma}(\bar{c}) \bigr)$ and
  Proposition~\ref{pro:moduleIso} shows that this agrees with
  multiplication in the deformed group ring.  Moreover, when
  $\ev_{2}^{*}(y^{b_{i}}) = 0$, we have $y^{c} * y^{b_{i}} = 0 = y^{c}
  \cdot y^{b_{i}}$.
  
  Next, consider the product $y^{v_{1}} * y^{v_{2}}$ where $v_{1},
  v_{2} \in \BBox(\bm{\Sigma})$.  If $\bar{v}_{1}$ and $\bar{v}_{2}$
  are not contained in a cone, then Remark~\ref{rem:nonempty} implies
  that $y^{v_{1}} * y^{v_{2}} = 0$ and \eqref{eq:definingProd} implies
  that $y^{v_{1}} \cdot y^{v_{2}} = 0$.  On the other hand, suppose
  the cone $\sigma \in \Sigma$ contains $\bar{v}_{1}$ and
  $\bar{v}_{2}$.  Let $v_{3} \in \BBox(\bm{\Sigma})$ be the element
  such that $\bar{v}_{3} \in \sigma(\bar{v}_{1}, \bar{v}_{2})$ and
  $v_{1} + v_{2} + v_{3} \equiv 0$; in other words, there exists
  integers $m_{i}$ such that $v_{1} + v_{2} + v_{3} = \sum_{\rho_{i}
    \in \sigma(\bar{v}_{1}, \bar{v}_{2}, \bar{v}_{3})} m_{i} b_{i}$
  and $1 \leq m_{i} \leq 2$.  Proposition~\ref{pro:K_0,3} shows that
  the component $\mathcal{K}_{v_{1}, v_{2}, v_{3}}$ is isomorphic to
  $\mathcal{X} \bigl( \bm{\Sigma / \sigma}(\bar{v}_{1}, \bar{v}_{2},
  \bar{v}_{3}) \bigr)$ and the evaluation map $\ev_{i}$ correspond to
  the closed embedding $\mathcal{X} \bigl( \bm{\Sigma /
    \sigma}(\bar{v}_{1}, \bar{v}_{2}, \bar{v}_{3}) \bigr) \rightarrow
  \mathcal{X} \bigl( \bm{\Sigma / \sigma}(\bar{v}_{i}) \bigr)$.  If
  $I$ is the set of indices $i$ such that $m_{i} = 2$, then
  Proposition~\ref{pro:virtual} shows that the virtual fundamental
  class on $\mathcal{X} \bigl( \bm{\Sigma / \sigma}(\bar{v}_{1},
  \bar{v}_{2}, \bar{v}_{3}) \bigr)$ is the product of the pullbacks of
  the divisor classes $y^{b_{i}}$ where $i \in I$.  Because of the
  degree shift, the class $y^{v_{i}} \in A_{orb}^{*} \bigl(
  \mathcal{X} (\bm{\Sigma}) \bigr)$ is identified with the class $1
  \in A^{*} \bigl( \mathcal{X} \bigl( \bm{\Sigma /
    \sigma}(\bar{v}_{i}) \bigr) \bigr)$ and $y^{v_{1}} * y^{v_{2}}$ is
  the image of the virtual fundamental class under the twisted
  evaluation map $\check{\ev}_{3}$.  In particular, if $J$ denotes the
  set of indices $i$ such that $\bar{b}_{i} \in \sigma(\bar{v}_{1},
  \bar{v}_{2})$ but $b_{i} \not\in \sigma(\bar{v}_{3})$, then
  unraveling the identification maps shows that $y^{v_{1}} *
  y^{v_{2}} = y^{\check{v}_{3}} \cdot \prod_{i \in I} y^{b_{i}} \cdot
  \prod_{j \in J} y^{b_{j}}$ where $\check{v}_{3}$ is the
  representation of $- v_{3}$ in $\in \BBox(\bm{\Sigma})$.  The factor
  $y^{\check{v}_{3}}$ arises from the involution $\iota \colon
  \mathcal{I} \bigl( \mathcal{X}(\bm{\Sigma}) \rightarrow \mathcal{I}
  \bigl( \mathcal{X}(\bm{\Sigma}) \bigr)$.  Since $\check{v}_{3} +
  \sum_{i \in I} b_{i} + \sum_{j \in J} b_{j} = v_{1} + v_{2}$, we
  conclude that $y^{v_{1}} * y^{v_{2}} = y^{v_{1}} \cdot y^{v_{2}}$.
\end{proof}

\section{Applications to Crepant Resolutions} \label{sec:crepant}
\setcounter{equation}{0}

In this section, we relate the orbifold Chow ring to the Chow ring of
a crepant resolution by showing that both rings are fibers of a flat
family.  This provide a new proof that the graded components of these
Chow rings have the same dimension.  On the other hand, we also
establish that these Chow rings are not generally isomorphic.

A rational simplicial fan $\Sigma$ with $n$ rays produces a canonical
stacky fan $\bm{\Sigma} := (N,\Sigma,\beta)$ where $N$ is the
distinguished lattice in the vector space containing $\Sigma$ and
$\beta \colon \mathbb{Z}^{n} \rightarrow N$ is the map defined by the
minimal lattice points on the rays.  Hence, there is a natural toric
Deligne-Mumford stack $\mathcal{X}(\bm{\Sigma})$ associated to every
simplicial toric variety $X(\Sigma)$.  Proposition~\ref{pro:coarse}
shows that $X(\Sigma)$ is the coarse moduli space of
$\mathcal{X}(\bm{\Sigma})$.

\begin{theorem} \label{thm:flat}
  Let $X(\Sigma)$ be a complete simplicial toric variety and let
  $\mathcal{X}(\bm{\Sigma})$ be the associated toric Deligne-Mumford
  stack.  If $\Sigma'$ is a regular subdivision of $\Sigma$ such that
  $X(\Sigma')$ is a crepant resolution of $X(\Sigma)$, then there is a
  flat family $T \rightarrow \mathbb{P}^{1}$ of schemes such that
  $T_{0} \cong \Spec A_{orb}^{*} \bigl( \mathcal{X}(\bm{\Sigma})
  \bigr)$ and $T_{\infty} \cong \Spec A^{*} \bigl( X(\Sigma') \bigr)$.
\end{theorem}

\begin{remark}
  Any regular subdivision $\Sigma'$ of $\Sigma$ induces a morphism
  $X(\Sigma') \rightarrow X(\Sigma)$; see Section~1.4 in
  \cite{Fulton}.  This morphism is a crepant resolution if and only if
  there is $\Sigma$-linear support function $h' \colon \mathbb{Q}^{d}
  \rightarrow \mathbb{Q}$ such that $h'(0) = 0$ and $h'(b_{i}) = -1$
  where $b_{i}$ for $1 \leq i \leq m$ are the minimal lattices points
  on the rays in $\Sigma'$; see Section~3.4 in \cite{Fulton}.
\end{remark}

\begin{proof}
  To establish this theorem, we construct a family of algebras over
  $\mathbb{P}^{1}$ such that the fiber over zero is isomorphic to
  $A_{orb}^{*} \bigl( \mathcal{X} (\bm{\Sigma}) \bigr)$ and the fiber
  over $\infty$ is isomorphic to $A_{\infty} \cong A^{*} \bigl(
  \mathcal{X} (\bm{\Sigma'}) \bigr)$.  We also prove that this family
  is flat outside of a Zariski closed subset of $\mathbb{P}^{1} - \{
  0, \infty \}$.  The required family $T \rightarrow \mathbb{P}^{1}$
  is obtained by extending our family over this finite set.
  
  We first introduce some notation.  Let $b_{1}, \dotsc, b_{n}$ be the
  minimal lattice points on the rays in $\Sigma$ and let $b_{n+1},
  \dotsc, b_{m}$ be the minimal lattice points on the additional rays
  in $\Sigma'$.  Since $X(\Sigma')$ is smooth, the lattice points
  $b_{1}, \dotsc, b_{m}$ generate the group $N$.  Hence, the ring
  $\mathbb{Q}[N]^{\bm{\Sigma}}$ is isomorphic to the quotient of the
  polynomial ring $S := \mathbb{Q}[y^{b_{1}}, \dotsc, y^{b_{m}}]$ by
  the binomial ideal $I_{2}$ which encodes the multiplication rules in
  \eqref{eq:definingProd}.  Fix a $\mathbb{Z}$\nobreakdash-basis
  $\theta_{1}, \dotsc, \theta_{d}$ for $N^{\star} :=
  \Hom_{\mathbb{Z}}(N, \mathbb{Z})$ and let $I_{1}$ be the ideal in
  $S$ generated by linear equations $\sum_{i=1}^{m} \theta_{j}(b_{i})
  \, y^{b_{i}}$ for $1 \leq j \leq d$.  Finally, the assumption that
  $\Sigma'$ is a regular subdivision of $\Sigma$ means that there is a
  $\Sigma'$\nobreakdash-linear support function $h \colon N
  \rightarrow \mathbb{Z}$ such that $h(b_{i}) = 0$ for $1 \leq i \leq
  n$, $h(b_{i}) > 0$ for $n+1 \leq i \leq m$ and $h(c_{1} + c_{2})
  \geq h(c_{1}) + h(c_{2})$ for all lattice points $c_{1}$, $c_{2}$
  lying in the same cone of $\Sigma$.  This inequality is strict
  unless $c_{1}$ and $c_{2}$ lie in the same cone of $\Sigma'$.
  
  To describe our family over $\mathbb{P}^{1} - \{ \infty \}$, let
  $\breve{I}_{1}$ be the ideal in $S[t_{1}]$ generated by
  $\sum_{i=1}^{m} \theta_{j}(b_{i}) \, y^{b_{i}} t_{1}^{h(b_{i})}$ for
  $1 \leq j \leq d$.  Since $h(b_{i}) = 0$ if and only if $1 \leq i
  \leq n$, Theorem~\ref{thm:main} implies that
  \[
  \frac{S[t_{1}]}{\breve{I}_{1} + I_{2} + \langle t_{1} \rangle} \cong
  \frac{S}{\langle \sum_{i=1}^{n} \theta_{j}(b_{i}) y^{b_{i}} : 1 \leq
    j \leq d \rangle + I_{2}} \cong A_{orb}^{*} \bigl( \mathcal{X}
  (\bm{\Sigma}) \bigr) \, .
  \]
  Moreover, it follows from the decomposition of $I_{2}$ in
  Proposition~4.8 in \cite{ES} that the sequence $\sum_{i=1}^{n}
  \theta_{j}(b_{i}) y^{b_{i}}$ for $1 \leq j \leq d$ forms a
  homogeneous system of parameters on $S/I_{2}$.  Lemma~4.6 in
  \cite{Stanley} shows that $S/I_{2}$ is a Cohen-Macaulay ring, so we
  deduce that $\sum_{i=1}^{n} \theta_{j}(b_{i}) y^{b_{i}}$ for $1 \leq
  j \leq d$ is a regular sequence.  Being a regular sequence is an
  open condition on the set of $d$-tuples of degree one elements in a
  finitely generated $\mathbb{Q}$-algebra.  Therefore, the Hilbert
  function of the family $S[t_{1}]/(\breve{I}_{1} + I_{2})$ over
  $\mathbb{Q}[t_{1}]$ is constant outside a finite set in
  $\mathbb{Q}^{*}$.
  
  For the family over $\mathbb{P}^{1} - \{ 0 \}$, let $\breve{I}_{2}$
  be the binomial ideal in $S[t_{2}]$ which encodes the product rule
  \[
  y^{c_{1}} \cdot y^{c_{2}} =
  \begin{cases}
    y^{c_{1} + c_{2}} \, t_{2}^{h(c_1+c_2)-h(c_1)-h(c_2)} & \text{if
      there exists $\sigma \in \Sigma$ such that $c_{1}, c_{2}
      \in \sigma$,}\\
    0 & \text{otherwise.}
  \end{cases}
  \] 
  Because $h(c_{1} + c_{2}) \geq h(c_{1}) + h(c_{2})$ and equality
  holds if and only if $c_{1}$ and $c_{2}$ lie in the same cone of
  $\Sigma'$, this product becomes
  \[
  y^{c_{1}} \cdot y^{c_{2}}=
  \begin{cases}
    y^{c_{1}+c_{2}} & \text{if there exists $\sigma' \in \Sigma'$ such
      that $c_{1},c_2 \in \sigma'$,}\\
    0 & \text{otherwise}
  \end{cases}
  \]
  over $t_{2} = 0$.  It follows that $S[t_{2}]/(\breve{I}_{2} +
  \langle t_{2} \rangle) \cong S/I_{\Sigma'}$, where $I_{\Sigma'}$ is
  the Stanley-Reisner ideal associated to $\Sigma'$, and
  Lemma~\ref{lem:basicChow} shows that $S[t_{2}]/( I_{1} +
  \breve{I}_{2} + \langle t_{2} \rangle) \cong A^{*} \bigl(
  \mathcal{X} (\bm{\Sigma}) \bigr)$.  As $\mathbb{Q}$-vector spaces,
  both $S/I_{2}$ and $S/I_{\Sigma'}$ have a basis consisting of the
  monomials in $S$ corresponding to lattice points in $N$.
  Proposition~4.8 in \cite{ES} implies that the sequence
  $\sum_{i=1}^{n} \theta_{j}(b_{i}) y^{b_{i}}$ for $1 \leq j \leq d$
  forms a homogeneous system of parameters on $S/I_{2}$, and
  Theorem~5.1.16 in \cite{BH} shows that this sequence is also a
  homogeneous system of parameters on $S/I_{\Sigma'}$.  Thus,
  $\sum_{i=1}^{m} \theta_{j}(b_{i}) y^{b_{i}}$ for $1 \leq j \leq d$
  is a regular sequence on both $S/I_{2}$ and $S/I_{\Sigma'}$ and the
  Hilbert functions of $S/(I_{1} + I_{2})$ and $S/(I_{1} +
  I_{\Sigma'})$ are equal.
  
  We combine the two one-parameter families by using the automorphisms
  $\varphi_{k}$ for $k = 1$, $2$ of $S[t_{k}^{\,}, t_{k}^{-1}]$
  defined by $\varphi_{k}(y^{b_{i}}) = y^{b_{i}} t_{k}^{h(b_{i})}$.
  Since $\varphi_{k}$ takes $I_{k} \cdot S[t_{k},t_{k}^{-1}]$ to
  $\breve{I}_{k} \cdot S[t_{k},t_{k}^{-1}]$, these automorphisms
  induce the following isomorphisms:
  \[
  \frac{S[t_{1}^{\,}, t_{1}^{-1}]}{\breve{I}_{1} + I_{2}} \cong
  \frac{S}{I_{1} + I_{2}} \cong \frac{S[t_{2}^{\,}, t_{2}^{-1}]}{I_{1}
    + \breve{I}_{2}} \, .
  \]
  Since a family of affine cones is a flat family if and only if the
  Hilbert function is constant (see Proposition~III-56 in \cite{EH}),
  we conclude that our family is flat on a Zariski open subset of
  $\mathbb{P}^{1}$ which contains both $0$ and $\infty$.
\end{proof}

\begin{remark}
  In analogy with Theorem~15.17 in \cite{Eisenbud}, the flat family
  constructed in the proof of Theorem~\ref{thm:flat} can be
  interpreted as a pair of Gr\"{o}bner deformations with respect to
  the appropriate weight orders connecting the ideal $\breve{I}_{1} +
  \breve{I}_{2}$ with its initial ideals $\langle \sum_{i=1}^{n}
  \theta_{j}(b_{i}) y^{b_{i}} : 1 \leq j \leq d \rangle +
  \breve{I}_{2}$ and $\breve{I}_{1} + I_{\Sigma'}$.
\end{remark}

We end with an example in which $A^{*} \bigl( X(\Sigma') \bigr)$ is
not isomorphic to $A_{orb}^{*} \bigl( \mathcal{X} (\bm{\Sigma})
\bigr)$.

\begin{example}
  Let $N = \mathbb{Z}^{2}$ and let $\Sigma \subseteq \mathbb{R}^{2}$
  be the complete fan in which the rays are generated by the lattice
  points $b_{1} := (1,0)$, $b_{2} := (0,-1)$ and $b_{3} := (-1,2)$.
  Hence, the toric variety $X(\Sigma)$ is the weighted projective
  space $\mathbb{P}(1,2,1)$ and the associated toric Deligne-Mumford
  stack is the quotient $[(\mathbb{C}^{3} - \{ 0 \})/ \mathbb{C}^{*}]$
  where the action is given by $(z_{1}, z_{2}, z_{3}) \cdot \lambda =
  (\lambda z_{1}, \lambda ^{2} z_{2}, \lambda z_{3})$.  If we simply
  write $x_{i}$ for the element $y^{b_{i}} \in \mathbb{Q}[N]$, then
  Theorem~\ref{thm:main} implies that
  \[
  A_{orb}^{*} \bigl( \mathcal{X}(\bm{\Sigma}) \bigr) \cong
  \frac{\mathbb{Q}[x_{1}, x_{2}, x_{3}, x_{4}]} {\bigl\langle x_{1}
    x_{3} - x_{4}^2, x_{2}x_{4}, x_{1} - x_{3}, - x_{2} + 2x_{3}
    \bigr\rangle} \cong \frac{\mathbb{Q}[x_{3}, x_{4}]} {\bigl\langle
    x_{3}^2 - x_{4}^{2}, x_{3}x_{4} \bigr\rangle} \, .
  \]
  Let $\Sigma'$ be the fan obtained from $\Sigma$ by inserting the ray
  generated by $b_{4} := (0,1)$.  It follows that $X(\Sigma')$ is the
  Hirzeburch surface $\mathbb{F}_{2}$, $X(\Sigma') \rightarrow
  X(\Sigma)$ is a crepant resolution (it blows down the $(-2)$-curve
  in $\mathbb{F}_{2}$), and Lemma~\ref{lem:basicChow} gives:
  \begin{align*}
  A^{*} \bigl( X(\Sigma') \bigr) &\cong \frac{\mathbb{Q}[x_{1}, x_{2},
    x_{3}, x_{4}]} {\bigl\langle x_{1}x_{3}, x_{2}x_{4}, x_{1} -
    x_{3}, - x_{2} + 2 x_{3} + x_{4} \bigr\rangle} \\
  &\cong \frac{\mathbb{Q}[x_{3}, x_{4}]} {\bigl\langle x_{3}^{2},
    2x_{3}x_{4} + x_{4}^{2} \bigr\rangle} = \frac{\mathbb{Q}[x_{3},
    x_{4}]} {\bigl\langle x_{3}^{2}, (x_{3}+x_{4})^{2} \bigr\rangle}
  \, .
  \end{align*}
  Since there is a degree one element $x \in A^{*} \bigl( X(\Sigma')
  \bigr)$ such that $x^{2} = 0$ and $A_{orb}^{*} \bigl(
  \mathcal{X}(\bm{\Sigma}) \bigr)$ does not contain such an element,
  we conclude that $A_{orb}^{*} \bigl( \mathcal{X}(\bm{\Sigma}) \bigr)
  \not\cong A^{*} \bigl( X(\Sigma') \bigr)$.
\end{example}

\providecommand{\bysame}{\leavevmode\hbox to3em{\hrulefill}\thinspace}

\end{document}